\documentclass{amsart}
\usepackage{righttag}
\usepackage{epsf}
\def\hepsffile{\leavevmode\epsffile}

\theoremstyle{plain}
\newtheorem{thm}{Theorem}[subsection]
\newtheorem{cor}[thm]{Corollary}
\newtheorem{lem}[thm]{Lemma}

\newtheorem{prop}[thm]{Proposition}

\theoremstyle{remark}
\newtheorem{rem}{Remark}
\newtheorem{remarks}{Remarks}

\theoremstyle{definition}
\newtheorem{defin}[thm]{Definition}

\newtheorem{emf}[thm]{}

\def\St{\protect\operatorname{St}}

\def\pr{\protect\operatorname{pr}}

\def\Z{{\mathbb Z}}

\def\R{{\mathbb R}}

\def\1{\hbox{\rm\rlap {1}\hskip.03in{\textrm I}}}
\def\Bbbone{{\rm1\mathchoice{\kern-0.25em}{\kern-0.25em}
        {\kern-0.2em}{\kern-0.2em}I}}

\begin{document} 
\font\eightrm=cmr8
\font\tenrm=cmr10
\font\eightit=cmti8
\centerline{{\em This article has been published in:\/}}
\centerline{{\em J. Knot Theory Ramifications {\bf 8} (1999), no. 1, pp.
71-97}}
\title[Arnold-type Invariants of Curves on Surfaces]
{Arnold-type Invariants of Curves on Surfaces}
\author[V.~Tchernov]{Vladimir Tchernov}
\address{D-MATH, HG G66.4, ETH Zentrum, CH-8092 Z\"urich, Switzerland}
\email{chernov@math.ethz.ch}

\begin{abstract}

Recently V.~Arnold introduced Strangeness and $J^{\pm}$ invariants of
generic 
immersions of an oriented circle to $\R^2$. 
Here these invariants are generalized to the case of
generic immersions of an oriented circle to an arbitrary surface $F$.
We explicitly describe all the invariants 
satisfying axioms, which naturally generalize 
the axioms used by V.~Arnold.

\end{abstract}
 
\maketitle

By a surface we mean any smooth two-dimensional manifold, possibly with
boundary.

\section{Introduction}
Consider the space $\mathcal F$ of all curves (immersions of an oriented circle) 
on a surface $F$. 
We call a
curve generic, if its only multiple points are double points of
transversal self-intersection. Nongeneric curves form a
discriminant hypersurface in $\mathcal F$. There are three main strata of the
discriminant. 
They are formed by curves with a triple point, curves with a 
self-tangency point, at which the velocity vectors of the two branches 
are pointing to the
same direction (direct self-tangency) and curves with a self-tangency point,
at which the velocity vectors of the two branches are pointing 
to the opposite directions
(inverse self-tangency). The union of these strata is dense in the
discriminant. In~\cite{Arnold} V.~Arnold associated a sign to a 
generic crossing of each of these strata. He also introduced $\St$, $J^+$ and
$J^-$ invariants of generic curves on $\R^2$, 
which change by a constant under a
positive crossing of the triple point, direct self-tangency and inverse 
self-tangency strata, respectively, and do not change under crossings of the
other two strata. These invariants give a lower bound for the number of 
crossings of each part of the discriminant, which are 
necessary to transform one
generic curve on $\R^2$ to another.

We construct generalizations of these invariants to the case when $F$ is any
surface (not necessarily $\R^2$). The fact, that for most surfaces 
the fundamental group is nontrivial, allows us to 
subdivide each of the three strata of the discriminant into pieces. We show
that this 
subdivision is natural from the point of view of the singularity theory. 
We take an integer valued function $\psi$ on the set of pieces obtained from
one stratum, and try to construct an invariant which increases by $\psi(P)$
under
a positive crossing of $P$ and does not change under crossings of the other
two strata. In an obvious sense $\psi$ is a derivative of such an invariant
and the invariant is an integral of $\psi$.
We introduce a condition on $\psi$ which is necessary and sufficient for
existence of such an invariant. Any integrable, in the sense above,
function $\psi$ defines this kind of an invariant up to an additive
constant.

If the surface $F$ is orientable, then the condition which corresponds to
the generalizations of $J^+$ and $J^-$ is automatically satisfied
and such an invariant
exists for any function $\psi$. For the generalization of $\St$ the
condition is not trivial. We reduce it to a
simple condition on $\psi$ which is sufficient for existence of such an
invariant. All these conditions are satisfied in the case of orientation
reversing curves.

When this work was complete and the main results of it were published as
preprints of Uppsala University~\cite{Tchernov1} and~\cite{Tchernov2} I
received a preprint of A.~Inshakov~\cite{Inshakov} 
containing similar results, obtained by him independently.

\section{Arnold's Invariants}
\subsection{Basic facts and definitions.}

A {\em curve\/} is a smooth immersion of (an oriented circle) $S^1$ into 
a (smooth) surface $F$.

A generic curve has only ordinary double points of transversal
self-intersection. All nongeneric curves form in the space of all
curves a {\em discriminant hypersurface\/}, or for short, the {\em
discriminant\/}.

A self-tangency point of (an oriented) curve is called a point of a {\em 
direct self-tangency\/}, if the velocity vectors at this point 
have the same direction; otherwise it is called a point of an
{\em inverse self-tangency\/}.

A {\em coorientation\/} of a smooth hypersurface in a functional space is
a local choice of one of the two parts, 
separated by this hypersurface, in a neighborhood of any of its points.
This part is called {\em positive\/}. 

The coorientation of the smooth part of a singular hypersurface is called
consistent, if the following consistency condition holds in a neighborhood of
any singular point of any stratum of codimension one on the hypersurface 
(of codimension two in the ambient functional space):

The intersection index of any generic small oriented closed loop with
a hypersurface (defined as a difference between the numbers of positive
and negative intersections) should vanish.

A hypersurface is called {\em cooriented\/}, if a consistent coorientation of
its smooth part is chosen, and {\em coorientable\/}, if such a coorientation
exists.

There are three parts of the discriminant hypersurface formed by the
curves having triple points, having direct self-tangencies, and having
inverse self-tangencies, respectively.

\begin{lem}[Arnold~\cite{Arnold}]\label{lemcoorient}
Each of these three parts of the discriminant hypersurface is coorientable.
\end{lem}

Consider a transversal crossing of the triple point stratum of the
discriminant.
A {\em vanishing triangle\/} is the triangle formed by the three branches
of the curve, corresponding to a subcritical or to a supercritical value of
the parameter near the triple point of the critical curve.

The {\em sign of a vanishing triangle\/} is defined by the following
construction. The orientation of the immersed circle defines the cyclic
order on the sides of the vanishing triangle. (It is the order of the
visits of the triple point by the three branches.) Hence, the sides of the
triangle acquire orientations induced by the ordering. But each side has
also its own orientation, which may coincide, or not, with the
orientations defined by the ordering.

For each vanishing triangle we define a quantity $q\in\{0,1,2,3\}$ to be the
number of sides of the vanishing triangle equally oriented by the ordering
and their direction. The {\em sign of the vanishing triangle\/} is $(-1)^q$.

\begin{defin}[of the sign of a crossing of a stratum]\label{coorient}
A transversal crossing of the direct self-tangency or of the inverse
self-tangency
stratum of the discriminant 
is {\em positive\/}, if the number of double points increases (by two).

A transversal crossing of the triple point stratum of the discriminant 
is {\em positive\/}, if the new-born
vanishing triangle is positive.
\end{defin}

\subsection{Invariants $\St, J^+$ and $J^-$.}
The {\em index} of an immersion of an oriented circle into an oriented plane
is the number of turns of the velocity vector. (The degree of the mapping
sending a point of the circle to the direction of the derivative of the
immersion at this point.)
The Whitney Theorem~\cite{Whitney} says that the connected components 
of the space of oriented planar curves are counted by the
indices of the curves.

Consider one of these components, that is, the space of immersions of a fixed
index.

\begin{thm}[Arnold~\cite{Arnold}]\label{St}
There exists a unique (up to an additive constant) invariant of generic
planar curves of a fixed index, whose value remains unchanged under crossings
of the self-tangency strata of the discriminant, 
but increases by one under the
positive crossing of the triple point stratum of the discriminant.
\end{thm}

\begin{figure}[htbp]
 \begin{center}
  \epsfxsize 10cm
  \hepsffile{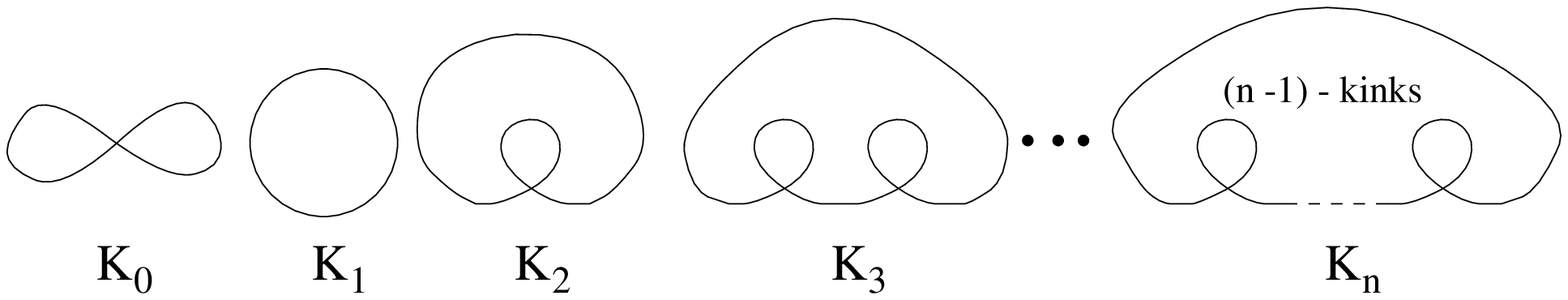}
 \end{center}
\caption{}\label{simple.fig}
\end{figure}

This invariant is denoted by $\St$ (from Strangeness), 
when normalized
by the following conditions:
\begin{equation}
\St(K_0)=0, \text{   }\St(K_{i+1})=i\text{   }(i=0,1,\dots),
\end{equation}
where $K_0$ is the figure eight curve and $K_{i+1}$ is the simplest curve
with $i$ double points (see~Figure~\ref{simple.fig}). The curve $K_j$ has index $\pm
j$, depending on the orientation. 

\begin{thm}[Arnold~\cite{Arnold}]\label{Jplus}
There exists a unique (up to an additive constant) invariant of generic
planar curves of a fixed index, whose value remains unchanged under a crossing of
the inverse self-tangency or of the triple point strata of the discriminant, 
but increases by two under the
positive crossing of the direct self-tangency stratum of the discriminant.
\end{thm}

This invariant is denoted by $J^+$, when normalized
by the following conditions:
\begin{equation}
J^+(K_0)=0, \text{   }J^+(K_{i+1})=-2i\text{   }(i=0,1,\dots),
\end{equation}
where $K_0$ and $K_{i+1}$ are the curves shown in~Figure~\ref{simple.fig}.

\begin{thm}[Arnold~\cite{Arnold}]\label{Jminus}
There exists a unique (up to an additive constant) invariant of generic
planar curves of a fixed index, whose value remains unchanged under a 
crossing of the direct self-tangency or of the triple point strata of the
discriminant, 
but decreases by two under the
positive crossing of the inverse self-tangency stratum of the discriminant.
\end{thm}

This invariant is denoted by $J^-$, when normalized
by the following conditions:
\begin{equation}
J^-(K_0)=-1, \text{   }J^-(K_{i+1})=-3i\text{   }(i=0,1,\dots),
\end{equation}
where $K_0$ and $K_{i+1}$ are the curves shown in~Figure~\ref{simple.fig}.

These normalizations of the three invariants were chosen~\cite{Arnold} 
to make them independent of the orientation of the parameterizing circle 
and additive under the connected summation of planar curves.

\section{Strangeness-type Invariant of Curves on Surfaces.}
\subsection{Natural decomposition of the triple point stratum.}
\begin{defin}\label{generictriple}
Let $F$ be a surface. 
We say that a curve $\xi\subset F$ with a triple point $q$ is a 
{\em generic curve with a triple point\/}, if its only nongeneric singularity is
this triple point, at which every two branches are transverse to each other.
\end{defin}
 
\begin{emf}
Let $F$ be a surface. Let $B_3$ be a bouquet of three oriented circles
with a fixed
cyclic order on the set of them, and let $b$ be the base point of $B_3$.
Let $s:S^1\rightarrow F$ be a generic curve with a triple point $q$.

Let $\alpha:S^1\rightarrow B_3$ be a continuous
mapping such that:

a) $\alpha (s^{-1}(q))=b$.

b) $\alpha$ is injective on the complement of $s^{-1}(q)$.

c) The orientation induced by $\alpha$ on $B_3\setminus b$
coincides with the orientation of the circles of $B_3$.

d) The cyclic order induced on the set of circles of $B_3$ by traversing
$\alpha(S^1)$ according to the orientation of $S^1$ 
coincides with the fixed one.

The mapping $\phi:B_3\rightarrow F$ such that $s=\phi\circ\alpha$ is called
an {\em associated\/} with $s$ mapping of $B_3$.

Note, that the free homotopy class of the mapping of $B_3$ to $F$ realized
by $\phi$ is well
defined, modulo an automorphism of $B_3$ which preserves the orientation of
the circles and
the cyclic order on the set of them.

\end{emf}

\begin{defin}[of $T$-equivalence]
Let $s_1$ and $s_2$ be two generic curves with a triple point
(see~\ref{generictriple}). 
We say, that these curves are
{\em $T$-equivalent\/}, if 
there exist associated with them mappings of $B_3$ which are
free homotopic. The triple point stratum is naturally
decomposed into parts corresponding to different 
$T$-equivalence classes. 

We denote by $[s]$ the $T$-equivalence class corresponding to $s$, 
a generic curve
with a triple point. 
We denote by $\mathcal T$ the set of all the $T$-equivalence classes.
\end{defin}

\subsection{Axiomatic description of $\overline{\St}$.}
A natural way to introduce $\St$ type invariants of generic curves
on a surface $F$ is to take a function $\psi:\mathcal T\rightarrow \Z$ and to
construct an invariant of generic curves from a fixed
connected component $\mathcal C$ of
the space $\mathcal F$ (of all the curves on $F$) such that:

1. It does not change under crossings of the self-tangency 
strata of the discriminant.

2. It increases by $\psi ([s])$
under a positive crossing of the part of the triple point stratum, which
corresponds to a $T$-equivalence class $[s]$.

If for a given function $\psi:\mathcal T \rightarrow \Z$
there exists such an invariant of curves
from $\mathcal C$, then
we say that there exists a $\overline{\St}$ invariant of
curves in $\mathcal C$, which is an {\em integral\/} of $\psi$.
Such $\psi$ is said to be $\overline{\St}$-{\em integrable\/} in $\mathcal C$.

\begin{emf}\label{obstruct}{\em Obstructions for the integrability.\/} 
Fix $\psi:\mathcal T\rightarrow \Z$.
Let $\xi\subset F$ be a generic curve and $\gamma\subset \mathcal C$ be a
generic loop starting at $\xi$. 
We denote by $I^{\xi}$ the set of moments, when $\gamma$
crosses the triple point stratum.  We denote by
$\{t_i^{\xi}\}_{l\in I^{\xi}}$ the
$T$-equivalence classes corresponding to the parts of the stratum,
where the crossings occur,
and by $\{\sigma_i^{\xi}\}_{i\in
I^{\xi}}$ the signs of the crossings. Put
\begin{equation}
\Delta_{\overline{\St}}(\gamma)
=\sum_{i\in I^{\xi}}\sigma_i^{\xi}\psi(t^{\xi}_i)
\end{equation}

We call
$\Delta_{\overline{\St}}(\gamma)$ the 
{\em change of $\overline{\St}$ along $\gamma$.\/}
If $\Delta_{\overline{\St}}(\gamma)=0$, then $\psi$ is said to be 
{\em integrable along\/} $\gamma$. 
It is clear, that if a function $\psi$ is integrable in $\mathcal C$, then it is
integrable along any generic loop $\gamma\subset \mathcal C$.

Below we describe the two loops $\gamma_1$ and $\gamma_2$ such that the
integrability of $\psi$ along them implies integrability of $\psi$ in $\mathcal
C$. In a sense, the changes along these loops are the only obstructions for
the integrability. 
(The loop $\gamma_2$ is going to be well defined (and needed) only in the
case of $F$ being a Klein bottle and $\mathcal C$ consisting of orientation
reversing curves on it.)
\end{emf}

\begin{emf}\label{gamma1}{\em Loop $\gamma_1$.\/}
Let $\xi\in C$ be a generic curve and $\gamma_1\subset \mathcal C$ be the loop
starting at $\xi$, which is constructed below.

Deform $\xi$ along a generic path $t$ in $\mathcal
F$ to
get two opposite kinks,  as it shown in Figure~\ref{twokink.fig}. Make
the first kink very small and slide it along the curve (in such a way that at
each moment of time points of $\xi$ located outside of a small neighborhood
of the kink do not move)  till it comes back.   (See
Figure~\ref{factorst.fig}.) 
Finally deform $\xi$ to its original shape along $t^{-1}$. 

Note, that if $\xi$ represents an orientation reversing loop on $F$, 
then the kink slides twice along $\xi$ before
it returns to the original position.

\begin{figure}[htbp]
 \begin{center}
  \epsfxsize 10cm
  \hepsffile{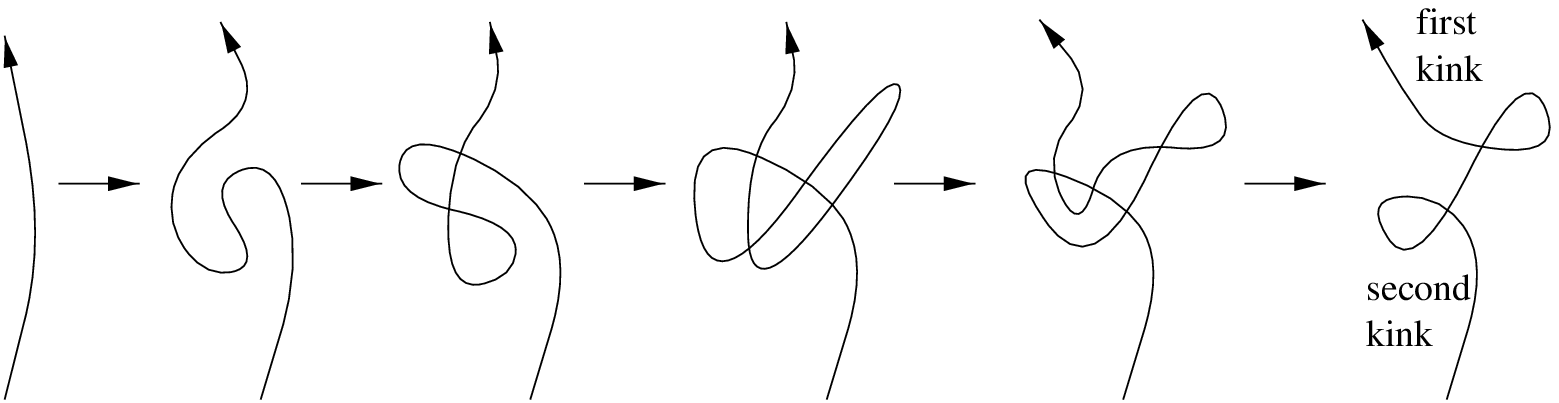}
 \end{center}
\caption{}\label{twokink.fig}
\end{figure}

\begin{figure}[htbp]
 \begin{center}
  \epsfxsize 8cm
  \hepsffile{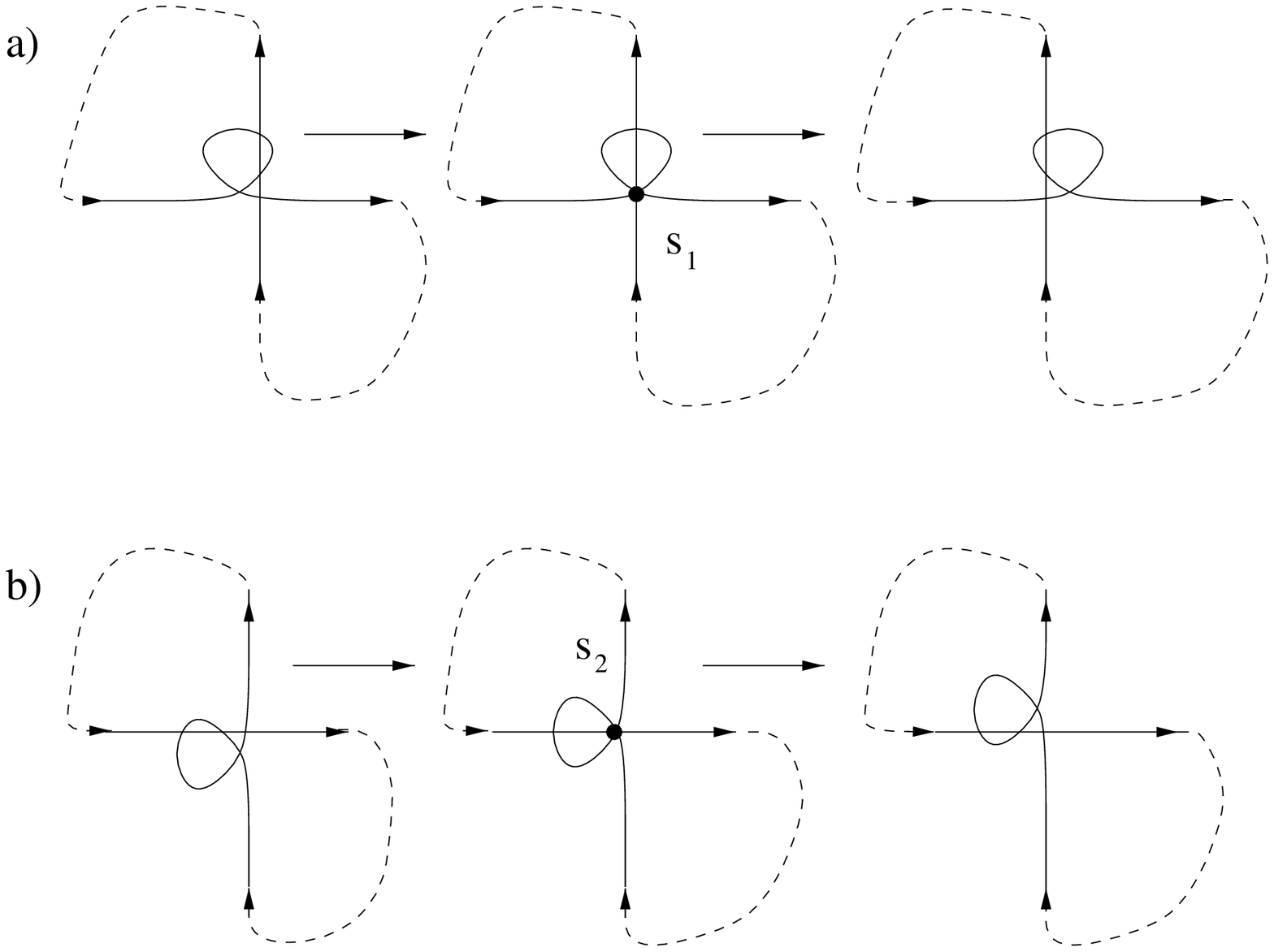}
 \end{center}
\caption{}\label{factorst.fig}
\end{figure}

\end{emf}

\begin{emf}\label{gamma2}{\em Loop $\gamma_2$.\/} 
Let $\xi$ be a generic orientation preserving 
curve on the Klein bottle $K$. Let $\gamma_2\subset \mathcal C$ be the loop
starting at $\xi$ which is constructed below.

Consider $K$ as a quotient of a rectangle modulo the identification on its
sides shown in Figure~\ref{klein.fig}. Let $p$ be the orientation
covering $T^2\rightarrow K$. There is a loop $\alpha$ in the space of all
autodiffeomorphisms of $T^2$,
which is the sliding of $T^2$ along the unit
vector field parallel to the lifting
of the curve $c\subset K$ 
(see Figure~\ref{klein.fig}). Since $\xi$ is an orientation preserving curve 
it can
be lifted to a curve $\xi '$ on $T^2$.
The loop $\gamma_2$ is the composition of $p$
and of the sliding of $\xi '$ induced by $\alpha$. (To make $\gamma_2$
well-defined for each $\xi$ we choose which one of the two possible 
liftings of $\xi$ to a curve on $T^2$ is $\xi'$.)

\begin{figure}[htbp]
 \begin{center}
  \epsfxsize 4cm
  \hepsffile{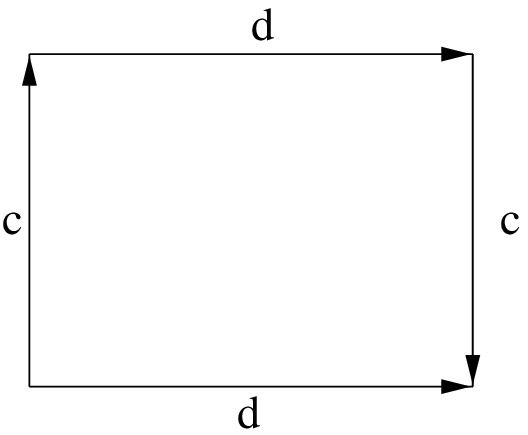}
 \end{center}
\caption{}\label{klein.fig}
\end{figure}

\end{emf}

\begin{thm}\label{barST}
Let $F$ be a surface (not necessarily compact or orientable), 
$\mathcal T$ be the set of all the $T$-equivalence
classes, $\mathcal C$ be a connected component of $\mathcal F$ and $\xi\in \mathcal C$
be a generic curve. 
Let $\psi:\mathcal T\rightarrow \Z$ be a function. 

Then the following two statements \textrm{I} and \textrm{II} are equivalent.

\begin{description}
\item[\textrm{I}] There exists an invariant $\overline{\St}$
of generic curves from $\mathcal C$ which is an integral of $\psi$.

\item[\textrm{II}] If $F\neq K$ (Klein bottle) then $\psi$ is integrable along 
the loop $\gamma_1\subset \mathcal C$ starting at $\xi$.

If $F=K$ and $\mathcal C$ consists of orientation
reversing curves on $K$, then $\psi$ is integrable along the loop
$\gamma_1\subset \mathcal C$ starting at $\xi$. 
If $F=K$ and $\mathcal C$ consists of orientation
preserving curves on $K$, then $\psi$ is integrable along the loops
$\gamma_1, \gamma_2\subset \mathcal C$ starting at $\xi$. 
\end{description}
\end{thm}

For the Proof of Theorem~\ref{barST} see
Section~\ref{pfbarST}.

\begin{remarks}  
If for a given function $\psi:\mathcal T\rightarrow \Z$ there exists
an invariant $\overline{\St}$ which is an integral of $\psi$, 
then it is unique up to an additive
constant. (This statement follows from the proof of Theorem~\ref{barST}.) 

Note, that 
if statement \textrm{II} holds for one
generic $\xi\in \mathcal C$, then statement \textrm{I} holds, which implies that 
\textrm{II} holds for all generic $\xi' \in\mathcal C$. 

A straightforward modification of the proof of Theorem~\ref{barST} shows
that it holds for $\psi$ 
taking values in any torsion free Abelian group.

The connected components of $\mathcal F$ admit a rather simple description.
One can show (cf.~\ref{H-principle}) that they 
are naturally identified with the connected 
components of the space of free loops in $STF$ (the spherical tangent bundle
of the surface) or, which is the same, with the conjugacy classes of
$\pi_1(STF).$
\end{remarks}

\begin{emf}{\em Cases, when $\psi$ is automatically integrable.\/}
Theorem~\ref{barST} says, that in the cases of orientable $F$, 
or of $\mathcal C$ consisting of orientation reversing curves on $F$,
integrability of $\psi$ along $\gamma_1$ is
sufficient for the existence of the $\overline{\St}$ invariant.

Clearly, all the crossings of the triple point stratum, which occur along
$\gamma_1$  (sliding of a kink along $\xi$) 
happen, when the kink passes through a double
point of $\xi$. (See Figure~\ref{factorst.fig}.) 

If $F$ is orientable, then the kink passes twice 
through each double point. A straightforward check shows that the signs of 
the corresponding triple point stratum 
crossings are opposite. The mappings of $B_3$ associated with these
crossings are different by an orientation preserving automorphism of
$B_3$, which does not preserve the cyclic order on the circles. 
For both crossings the restriction of an associated mapping 
to one of the circles of $B_3$ represents a contractible loop.
Thus, if $F$ is orientable then 
statement {\textrm II} (and hence statement {\textrm I}) 
of Theorem~\ref{barST} is true for any function $\psi$,
provided that it takes the same value on any two $T$-equivalence classes, 
for which there exist mappings $\phi _1, \phi _2$ of $B_3$ representing them, 
such that: 

a) The restriction of $\phi _i$ (${i\in \{1,2\}}$) to one of the circles of $B_3$
represents a contractible loop on $F$.

b) There exists $\alpha$, an orientation preserving 
automorphism of $B_3$ (not preserving the cyclic order on the
circles) such that $\phi _1=\phi _2\circ\alpha$.

If $\xi$ represents an orientation reversing loop on
$F$, then the kink has to slide twice along $\xi$ before it comes to
the original position. 
Thus, it passes four times through each double point
of $\xi$. One can show, that the corresponding crossings of the triple
point stratum can
be subdivided into two pairs, such that the $T$-equivalence classes
corresponding to the crossings inside each pair are equal and the signs
of the two crossings in each pair are opposite. Thus, if $\xi$ represents an
orientation reversing loop on $F$, then statement {\textrm II} (and hence
statement {\textrm I}) of Theorem~\ref{barST} is true for
any function $\psi:\mathcal T\rightarrow \Z$. 
Another way of proving this is
based on the fact that for such $\xi$ the loop 
$\gamma_1=1\in \pi_1(\mathcal F, \xi)$, see
Section~\ref{H-principle}.   
\end{emf}

\begin{rem} The following example shows that on nonorientable surfaces even
a constant (nonzero) function $\psi$ is not necessarily integrable.

Let $F$ be a nonorientable surface. Let $\xi\in\mathcal C$ 
be an orientation
preserving curve with a single double point $x$ which separates $\xi$ into two 
orientation reversing loops. Then $\Delta_{\overline{\St}}(\gamma_1)\neq 0$
for a constant nonzero function $\psi:\mathcal T\rightarrow \Z$. (The reason is
that the signs of the two crossings of the triple point stratum 
corresponding to a kink passing through $x$ have the same sign.)
\end{rem}

\begin{emf}\label{splitST}{\em Connection with the standard $\St$-invariant.\/}
Since $\R^2$ is simply connected, 
there is just one $T$-equivalence class of
singular curves on $\R^2$.
Thus, the construction of $\overline{\St}$ does not give anything new in the
classical case of planar curves. 
\end{emf}

\subsection{Singularity theory interpretation of $\overline{\St}$ for
orientable $F$.}

\begin{defin}[of $\overline T$-equivalence]
Let $S^1(3)$ be the configuration space of unordered triples of distinct 
points on
$S^1$. Consider a space $S^1(3)\times\mathcal F$. Let $\mathcal M$ be the subspace of 
$S^1(3)\times\mathcal F$ consisting of $t\times f\in  S^1(3)\times\mathcal F$, such
that $f$ maps the three points from $t$ to one point on $F$. (This is a sort
of singularity resolution for strata involving points of multiplicity 
greater than two.)

We say that $m_1,m_2\in \mathcal M$ are $\overline T$-equivalent, 
if they belong to the same path connected component of $\mathcal M$. 
For $s$,
a generic curve with a triple point (see~\ref{generictriple}), 
there is a
unique $\overline T$-equivalence class associated to it. We denote this
class by $[\bar s]$. 

Thus, the 
$\overline T$-equivalence relation induces a decomposition  of 
the triple point stratum of the discriminant hypersurface.
\end{defin}

Let $\overline{\mathcal T}$ be the set of all the $\overline T$-equivalence classes. 
There is a
natural mapping $\phi:\overline{\mathcal T}\rightarrow \mathcal T$.  
It maps $\bar t\in \overline{\mathcal T}$ to such $t\in \mathcal T$, that there 
exists $s$, a generic curve with a triple point, 
for which $[s]=t$ and 
$[\bar s]=\bar t$. 

Let $F$ be an orientable surface.
Let $\mathcal C$ be a connected component of $\mathcal F$
and $\overline {\mathcal T_{\mathcal C}}\subset \overline{\mathcal T}$ be the set of
all the $\overline
T$-equivalence classes, corresponding to generic 
curves (from $\mathcal C$) with a triple point.

\begin{thm}\label{discrST}
The mapping $\phi\big|_{\overline{\mathcal T_{\mathcal C}}}$ is injective.
\end{thm} 

For the Proof of Theorem~\ref{discrST} see Section~\ref{pfdiscrST}.

\begin{emf}{\em Interpretation of $\overline{\St}$.\/}
Let $F$ be an orientable surface.
Let $\mathcal C$ be a connected component of $\mathcal F$. 
Let $\tilde {\St}$ be an invariant of generic curves from $\mathcal C$ such
that:
 
a) It does not change under crossings of the self-tangency strata of the
discriminant.  

b) Under the positive crossing of a part of the triple point stratum 
of the discriminant it increases by a constant depending only 
on the $\overline T$-equivalence class corresponding to this part of the 
stratum. 

Theorem~\ref{discrST} implies, that this $\tilde {\St}$ invariant is an
$\overline{\St}$ invariant for some choice of the function $\psi:\mathcal
T\rightarrow \Z$.
\end{emf}

\section{$J^+$-type Invariant of Curves on Surfaces.}
\subsection{Natural decomposition of the direct self-tangency point stratum.}

\begin{defin}\label{genericdirect}
Let $F$ be a surface. 
We say that a curve $\xi\subset F$ with a direct self-tangency point
$q$ is a 
{\em generic curve with a direct self-tangency point\/}, 
if its only nongeneric singularity is
this point.
\end{defin}

\begin{emf}
Let $F$ be a surface. 
Let $B_2$ be a bouquet of two oriented circles, and $b$ be its base point. 
Let $s:S^1\rightarrow F$ be a generic curve with a direct self-tangency point 
$q$. It can be lifted to the mapping $\bar s$
from the oriented circle to $STF$ (the spherical tangent bundle of F),
which sends a point $p\in S^1$ to the point in $STF$ corresponding to the
direction of the velocity vector of $s$ at $s(p)$. 
(Note, that $q$ lifts to a double point
$\bar q$ of
$\bar s$.) 

Let $\alpha:S^1\rightarrow B_2$ be a continuous mapping such that:

a) $\alpha (\bar s^{-1}(\bar q))=b$.

b) $\alpha$ is injective on the complement of $\bar s^{-1}(\bar q)$.

c) The orientation induced by $\alpha$ on $B_2\setminus b$ coincides with
the orientation of the circles of $B_2$.

The mapping $\phi:B_2\rightarrow STF$ such that 
$\bar s=\phi\circ\alpha$ is called
an {\em associated\/} with $s$ mapping of $B_2$.

Note, that the free homotopy class of a mapping of $B_2$ to $STF$ realized by 
$\phi$ is well defined, modulo the orientation preserving automorphism of
$B_2$ which interchanges the circles.
\end{emf}

\begin{defin}[of $T^+$-equivalence]
Let $s_1$ and $s_2$ be two generic curves with a point of direct
self-tangency (see~\ref{genericdirect}). 
We say that these curves are
{\em $T^+$-equivalent\/} if there exist associated with the two of them
mappings of $B_2$, which are free homotopic. 
The direct self-tangency point stratum is naturally
decomposed into parts corresponding to different 
$T^+$-equivalence classes. 

We denote by $[s^+]$ the $T^+$-equivalence class corresponding to 
$s$, a generic curve with a point of direct self-tangency.
We denote by $\mathcal T^+$ the set of all the $T^+$-equivalence classes.

\end{defin}

\subsection{Axiomatic description of $\overline{J^+}$.}
A natural way to introduce $J^+$ type invariant of generic curves
on a surface $F$ is to take a function $\psi:\mathcal T^+\rightarrow \Z$ and to
construct an invariant of generic curves from a fixed
connected component $\mathcal C$ of
the space $\mathcal F$ (of all the curves on $F$) such that:

1. It does not change under crossings of the triple point or of the inverse
self-tangency strata of the discriminant.

2. It increases by $\psi ([s^+])$
under a positive crossing of the part of the direct self-tangency point 
stratum, which corresponds to a $T^+$-equivalence class $[s^+]$.

If for a given function $\psi:\mathcal T^+ \rightarrow \Z$
there exists such an invariant of curves
from $\mathcal C$, then
we say that there exists a $\overline{J+}$ invariant of
curves in $\mathcal C$, which is an {\em integral\/} of this function.
Such $\psi$ is said to be $\overline{J^+}$-{\em integrable\/} in $\mathcal C$.

Similarly to~\ref{obstruct} 
we introduce the notions of the {\em change of $J^+$ along a loop\/} and of
the {\em integrability of $\psi$ along a loop.\/}

\begin{thm}\label{barJ+}
Let $F$ be a surface (not necessarily compact or orientable), 
$\mathcal T^+$ be the set of all the $T^+$-equivalence
classes, $\mathcal C$ be a connected component of $\mathcal F$ and $\xi\in \mathcal C$
be a generic curve. 
Let $\psi:\mathcal T^+\rightarrow \Z$ be a function. 

Then the following two statements \textrm{I} and \textrm{II} are equivalent.

\begin{description}
\item[\textrm{I}] There exists an invariant $\overline{J^+}$
of generic curves from $\mathcal C$, which is an integral of $\psi$.

\item[\textrm{II}] If $F\neq K$ (Klein bottle) then $\psi$ is integrable along 
the loop $\gamma_1\subset \mathcal C$ starting at $\xi$.

If $F=K$ and $\mathcal C$ consists of orientation
reversing curves on $K$, then $\psi$ is integrable along the loop
$\gamma_1\subset \mathcal C$ starting at $\xi$. 
If $F=K$ and $\mathcal C$ consists of orientation
preserving curves on $K$, then $\psi$ is integrable along the loops
$\gamma_1, \gamma_2\subset \mathcal C$ starting at $\xi$. 
\end{description}
\end{thm}

For the Proof of Theorem~\ref{barJ+} see
Section~\ref{pfbarJ+}.

\begin{remarks}  
If for a given function $\psi:\mathcal T^+\rightarrow \Z$ there exists
an invariant $\overline{J^+}$ which is an integral of $\psi$, 
then it is unique up to an additive
constant. (This statement follows from the proof of Theorem~\ref{barJ+}.) 

Note, that 
if statement \textrm{II} holds for one
generic $\xi\in \mathcal C$, then statement \textrm{I} holds, which implies that 
\textrm{II} holds for all generic $\xi' \in\mathcal C$.

A straightforward modification of the proof of Theorem~\ref{barJ+} shows
that it holds for $\psi$ 
taking values in any torsion free Abelian group.

The connected components of $\mathcal F$ admit a rather simple description.
One can show (cf.~\ref{H-principle}) that they 
are naturally identified with the connected 
components of the space of free loops in $STF$ (the spherical tangent bundle
of the surface) or, which is the same, with the conjugacy classes of
$\pi_1(STF).$
\end{remarks}

\begin{emf}\label{trivialJ+}{\em Cases, when $\psi$ is automatically
integrable.\/}
Theorem~\ref{barJ+} says, that in the cases of orientable $F$, 
or of $\mathcal C$ consisting of orientation reversing curves on $F$,
integrability of $\psi$ along $\gamma_1$ is
sufficient for the existence of the $\overline{J^+}$ invariant.

Clearly, all the crossings of the direct self-tangency 
stratum, which occur along $\gamma_1$ (sliding of a kink along $\xi$) 
happen, when the kink passes through a double
point of $\xi$. 

If $F$ is orientable, then the kink passes twice 
through each double point of
$\xi$. A straightforward check shows, that the signs of the corresponding
direct self-tangency stratum 
crossings are opposite, and the $T^+$-equivalence classes corresponding to
them are equal. 
Thus, if $F$ is orientable, then statement {\textrm II} (and hence statement
{\textrm I})
of Theorem~\ref{barJ+} is true for any function $\psi:\mathcal T^+\rightarrow \Z$.

If $\xi$ is an orientation reversing curve on
$F$, then the kink slides twice along $\xi$, before it comes to
its original position. Thus, it passes four times through each double point
of $\xi$. One can show, that the corresponding four crossings of the
direct self-tangency stratum can
be subdivided into two pairs, such that the $T^+$-equivalence classes
corresponding to the crossings inside the same pair are equal and the signs
of the two crossings in each pair are opposite. Thus, if $\xi$ represents an
orientation reversing loop on $F$, then statement {\textrm II} (and hence
statement {\textrm I}) of
Theorem~\ref{barJ+} is true for
any function $\psi:\mathcal T^+\rightarrow \Z$. Another way of proving this is
based on the fact, that for such $\xi$ the loop 
$\gamma_1=1\in \pi_1(\mathcal F, \xi)$, see
Section~\ref{H-principle}. 
\end{emf}

\begin{rem} Similarly to the case of $\overline{\St}$, 
even a constant (nonzero) function $\psi$ 
is not necessarily integrable in the
case of orientation preserving $\xi$ on a nonorientable surface $F$.  
\end{rem}

\begin{emf}\label{splitJ+}{\em Connection with the standard
$J^+$-invariant.\/}
Since $\pi_1(ST\R^2)=\Z$, there are countably many $T^+$-equivalence 
classes of
singular curves on $\R^2$, 
which can be obtained from a curve of the fixed index.
(Note, that the index of a curve $\xi$ defines the connected component
of the space of all curves on $\R^2$, which $\xi$ belongs to.)
Thus, the construction of $\overline{J^+}$ gives rise to a splitting of the
standard $J^+$ invariant of V.~Arnold. This is the splitting introduced by
V.~Arnold~\cite{Arnoldsplit} in the case of planar curves of index zero 
and generalized to the case of arbitrary planar curves by F.~Aicardi~\cite{Aicardi}.
\end{emf}

\subsection{Singularity theory interpretation of $\overline{J^+}$ for
orientable $F$.}

\begin{defin}[of $\overline{T^+}$-equivalence]
Let $S^1(2)$ be the configuration space of unordered pairs of distinct 
points on
$S^1$. Consider a space $S^1(2)\times\mathcal F$. Let $\mathcal M^+$ be the subspace of 
$S^1(2)\times\mathcal F$ consisting of $t\times f\in S^1(2)\times\mathcal F$, such
that $f$ maps the two points from $t$ to one point on $F$ and the
velocity vectors of $f$ at these two points have the same direction.
(This is a sort of singularity resolution for the strata, involving points of
direct self-tangency.)

We say that $m_1^+$ and $m_2^+$ from $\mathcal M^+$ are $\overline{T^+}$-equivalent, 
if they belong to the same path connected component of $\mathcal M^+$. 

Clearly, for
$s$, a generic curve with a direct self-tangency point
(see~\ref{genericdirect}), 
there is a
unique $\overline{T^+}$-equivalence class associated with it. We denote this
class by $[\overline {s^+}]$. 
Thus, the 
$\overline{T^+}$-equivalence relation induces a decomposition of 
the direct self-tangency 
point stratum of the discriminant hypersurface.
\end{defin}

Let $\overline{\mathcal T^+}$ be the set of all the $\overline{T^+}$-equivalence classes. 
There is a
natural mapping $\phi:\overline{\mathcal T^+}\rightarrow \mathcal T^+$.  
It maps $\overline{t^+}\in \overline{\mathcal T^+}$ to such $t^+\in \mathcal T^+$, 
that there 
exists $s$ (a generic curve with a direct self-tangency point), 
for which $[s^+]=t^+$ and 
$[\overline {s^+}]=\overline{t^+}$. 

Let $F$ be an orientable surface.
Let $\mathcal C$ be a connected component of $\mathcal F$
and $\overline{\mathcal T^+_{\mathcal C}}\subset \overline{\mathcal T^+}$, be the set of
all the $\overline{T^+}$-equivalence 
classes corresponding to generic curves (from $\mathcal C$) 
with a point of direct self-tangency.

\begin{thm}\label{discrJ+}
The mapping $\phi\big|_{\overline{\mathcal T^+_{\mathcal C}}}$ is injective.
\end{thm} 

For the Proof of Theorem~\ref{discrJ+} see Section~\ref{pfdiscrJ+}.

\begin{emf}{\em Interpretation of $\overline{J^+}$.\/}
Let $F$ be an orientable surface.
Let $\mathcal C$ be a connected component of $\mathcal F$. 
Let $\tilde {J^+}$ be an invariant of generic curves from $\mathcal C$, such
that:
 
a) It does not change under crossings of the inverse self-tangency and
of the triple point strata of the
discriminant.  

b) Under the positive crossing of a part of the direct self-tangency point 
stratum 
of the discriminant it increases by a constant, depending only on the 
$\overline{T^+}$
equivalence class corresponding to this part of the stratum. 

Theorem~\ref{discrJ+} implies, that this $\tilde {J^+}$ invariant is a
$\overline{J^+}$ invariant for some choice of the function $\psi:\mathcal
T^+\rightarrow \Z$.
\end{emf}

\section{$J^-$-type Invariant of Curves on Surfaces.}
\subsection{Natural decomposition of the inverse self-tangency point stratum.}

\begin{defin}\label{genericinverse}
Let $F$ be a surface. 
We say that a curve $\xi\subset F$ with an inverse self-tangency point
$q$ is a 
{\em generic curve with an inverse self-tangency point\/}, 
if its only nongeneric singularity is
this point.
\end{defin}

\begin{emf}
Let $F$ be a surface. 
Let $B_2$ be a bouquet of two oriented circles, and $b$ be its base point. 
Let $s:S^1\rightarrow F$ be a generic curve with an inverse 
self-tangency point 
$q$. It can be lifted to the mapping $\bar s$
from the oriented circle to $PTF$ (the projectivized tangent bundle of F),
which sends a point $p\in S^1$ to the point in $PTF$ corresponding to the
tangent line containing the velocity vector of 
$s$ at $s(p)$. (Note, that $q$ lifts to a double point
$\bar q$ of
$\bar s$.) 

Let $\alpha:S^1\rightarrow B_2$ be a continuous mapping such that:

a) $\alpha (\bar s^{-1}(\bar q))=b$.

b) $\alpha$ is injective on the complement of $\bar s^{-1}(\bar q)$.

c) The orientation induced by $\alpha$ on $B_2\setminus b$ coincides with
the orientation of the circles of $B_2$.

The mapping $\phi:B_2\rightarrow PTF$ such that 
$\bar s=\phi\circ\alpha$ is called
an {\em associated\/} with $s$ mapping of $B_2$.

Note, that the free homotopy class of a mapping of $B_2$ to $PTF$ realized by 
$\phi$ is well defined, modulo the orientation preserving automorphism of
$B_2$ which interchanges the circles.
\end{emf}

\begin{defin}[of $T^-$-equivalence]
Let $s_1$ and $s_2$ be two generic curves with a point of an inverse
self-tangency (see~\ref{genericinverse}). 
We say, that these curves are
{\em $T^-$-equivalent\/}, if there exist associated with the two of
them mappings of $B_2$, which are free homotopic. 
The inverse self-tangency point stratum is naturally
decomposed into parts corresponding to different 
$T^-$-equivalence classes.

We denote by $[s^-]$ the $T^-$-equivalence class corresponding to $s$, a
generic curve with a point of an inverse self-tangency.
We denote by $\mathcal T^-$ the set of all the $T^-$-equivalence classes.

\end{defin} 

\subsection{Axiomatic description of $\overline{J^-}$.}

A natural way to introduce $J^-$ type invariant of generic curves
on a surface $F$ is to take a function $\psi:\mathcal T^-\rightarrow \Z$ and to
construct an invariant of generic curves from a fixed
connected component $\mathcal C$ of
the space $\mathcal F$ (of all the curves on $F$) such that:

1. It does not change under crossings of the triple point or of the direct
self-tangency strata of the discriminant.

2. It increases by $\psi ([s^-])$
under a positive crossing of the part of the inverse self-tangency point 
stratum, which corresponds to a $T^-$-equivalence class $[s^-]$.

If for a given function $\psi:\mathcal T^- \rightarrow \Z$
there exists such an invariant of curves
from $\mathcal C$, then
we say that there exists a $\overline{J^-}$ invariant of
curves in $\mathcal C$, which is an {\em integral\/} of this function.
Such $\psi$ is said to be $\overline{J^-}$-{\em integrable\/} in $\mathcal C$.

Similarly to~\ref{obstruct} 
we introduce the notions of the {\em change of $J^-$ along a loop\/} and of
the {\em integrability of $\psi$ along a loop.\/}

\begin{thm}\label{barJ-}
Let $F$ be a surface (not necessarily compact or orientable), 
$\mathcal T^-$ be the set of all the $T^-$-equivalence
classes, $\mathcal C$ be a connected component of $\mathcal F$ and $\xi\in \mathcal C$
be a generic curve. 
Let $\psi:\mathcal T^-\rightarrow \Z$ be a function. 

Then the following two statements \textrm{I} and \textrm{II} are equivalent.

\begin{description}
\item[\textrm{I}] There exists an invariant $\overline{J^-}$
of generic curves from $\mathcal C$, which is an integral of $\psi$.

\item[\textrm{II}] If $F\neq K$ (Klein bottle) 
then $\psi$ is integrable along 
the loop $\gamma_1\subset \mathcal C$ starting at $\xi$.

If $F=K$ and $\mathcal C$ consists of orientation
reversing curves on $K$, then $\psi$ is integrable along the loop
$\gamma_1\subset \mathcal C$ starting at $\xi$. 
If $F=K$ and $\mathcal C$ consists of orientation
preserving curves on $K$, then $\psi$ is integrable along the loops
$\gamma_1, \gamma_2\subset \mathcal C$ starting at $\xi$. 
\end{description}
\end{thm}

The Proof of Theorem~\ref{barJ-} is a straightforward generalization of the
Proof of Theorem~\ref{barJ+}. 

\begin{remarks}  
If for a given function $\psi:\mathcal T^-\rightarrow \Z$ there exists
an invariant $\overline{J^-}$ which is an integral of $\psi$, 
then it is unique up to an additive
constant. (This statement follows from the proof of Theorem~\ref{barJ-}.) 

Note, that 
if statement \textrm{II} holds for one
generic $\xi\in \mathcal C$, then statement \textrm{I} holds, which implies that 
\textrm{II} holds for all generic $\xi' \in\mathcal C$.

A straightforward modification of the proof of Theorem~\ref{barJ-} shows
that it holds for $\psi$ 
taking values in any torsion free Abelian group.

The connected components of $\mathcal F$ admit a rather simple description.
One can show (cf.~\ref{H-principle}) that they 
are naturally identified with the connected 
components of the space of free loops in $STF$ (the spherical tangent bundle
of the surface) or, which is the same, with the conjugacy classes of
$\pi_1(STF).$
\end{remarks}

\begin{emf}{\em Cases, when $\psi$ is automatically integrable.\/}
Similarly to~\ref{trivialJ+} one can
show, that statement {\textrm II} 
of Theorem~\ref{barJ-} is true for any function $\psi:\mathcal T^-\rightarrow
\Z$, provided that $F$ is orientable, or that $\mathcal
C$ consists of orientation reversing curves on $F$.  
\end{emf}

\begin{rem} Similarly to the case of $\overline{\St}$, 
even a constant (nonzero) function $\psi$ 
is not necessarily integrable in the
case of orientation preserving $\xi$ on a nonorientable surface $F$.  
\end{rem}

\begin{emf}\label{splitJ-}{\em Connection with the standard $J^-$-invariant.\/}
Since $\pi_1(PT\R^2)=\Z$, there are countably many $T^-$-equivalence classes of
singular curves, which can be obtained from a curve of the fixed index.
(Note, that the index of a curve $\xi$ defines the connected component
of the space of all curves on $\R^2$, which $\xi$ belongs to.)
Thus, the construction of $\overline{J^-}$ gives rise to a splitting of the
standard $J^-$ invariant of V.~Arnold. This splitting is analogous to the
splitting of $J^+$ 
introduced by
V.~Arnold~\cite{Arnoldsplit} in the case of planar curves of index zero 
and generalized to the case of arbitrary planar curves by F.~Aicardi~\cite{Aicardi}.
\end{emf}  

\subsection{Singularity theory interpretation of $\overline{J^-}$ for
orientable $F$.}

\begin{defin}[of $\overline{T^-}$-equivalence]
Let $S^1(2)$ be the configuration space of unordered pairs of distinct 
points on
$S^1$. Consider a space $S^1(2)\times\mathcal F$. Let $\mathcal M^-$ be the subspace of 
$S^1(2)\times\mathcal F$ consisting of $t\times f\in S^1(2)\times\mathcal F$, such
that $f$ maps the two points from $t$ to one point on $F$ and the
velocity vectors of $f$ at these two points have opposite directions.
(This is a sort of singularity resolution for the strata, involving points of
an inverse self-tangency.)

We say, that $m_1^-$ and $m_2^-$ from $\mathcal M^-$ are
$\overline{T^-}$-equivalent, 
if they belong to the same path connected component of $\mathcal M^-$. 

Clearly, for
$s$, a generic curve with a point of an inverse self-tangency
(see~\ref{genericinverse}), there is a 
unique $\overline{T^-}$-equivalence class associated to it. 
We denote it by
$[\overline{s^-}]$. 

Thus, the 
$\overline{T^-}$-equivalence relation induces a decomposition of 
the inverse self-tangency 
point stratum of the discriminant hypersurface.
\end{defin}

Let $\overline{\mathcal T^-}$ be the set of all the $\overline{T^-}$-equivalence classes.

There is a
natural mapping $\phi:\overline{\mathcal T^-}\rightarrow \mathcal T^-$.  
It maps $\overline{t^-}\in \overline{\mathcal T^-}$ to such $t^-\in \mathcal T^-$, 
that there 
exists $s$ (a generic curve with an inverse self-tangency point), 
for which $[s^-]=t^-$ and 
$[\overline{s^-}]=\overline{t^-}$. 

Let $F$ be an orientable surface.
Let $\mathcal C$ be a connected component of $\mathcal F$,
and $\overline{\mathcal T^-_{\mathcal C}}\subset \overline{\mathcal T^-}$ be the set of
all the $\overline {T^-}$-equivalence classes corresponding to generic 
curves (from $\mathcal C$) with a point of an inverse self-tangency.

\begin{thm}\label{discrJ-}
The mapping $\phi\big|_{\overline{\mathcal T^-_{\mathcal C}}}$ is injective.
\end{thm} 

The Proof of Theorem~\ref{discrJ-} is a straightforward generalization 
of the proof of Theorem~\ref{discrJ+}.

\begin{emf}{\em Interpretation of $\overline{J^-}$.\/}
Let $F$ be an orientable surface.
Let $\mathcal C$ be a connected component of $\mathcal F$. 
Let $\tilde {J^-}$ be an invariant of generic curves from $\mathcal C$, such
that:
 
a) It does not change under crossings of the direct self-tangency and
of the triple point strata of the
discriminant.  

b) Under the positive crossing of a part of the inverse self-tangency point 
stratum 
of the discriminant it increases 
by a constant, depending only on the $\overline{T^-}$
equivalence class corresponding to this part of the stratum. 

Theorem~\ref{discrJ-} implies, that this $\tilde {J^-}$ invariant is a
$\overline{J^-}$ invariant for some choice of the function $\psi:\mathcal
T^-\rightarrow \Z$.
\end{emf}

\section{Proofs}

\subsection{Proof of Theorem~\ref{barST}}\label{pfbarST}
Clearly, in order for $\overline{\St}$ to be well defined, 
the change of it along any generic closed loop in $\mathcal C$ should be
zero. Thus, we have proved that statement \textrm{I} implies statement \textrm{II}.

To prove that statement \textrm{II} implies statement \textrm{I} we imitate the
approach developed by V.~Arnold~\cite{Arnold} in the case of planar curves.

Fix any value of $\overline{\St}(\xi)\in
\Z$.
Let $\xi'\in \mathcal C$ be another generic curve. 
Take a generic path $p$ in $\mathcal C$, which connects 
$\xi$ with $\xi'$. When we go along this path we see a sequence of
crossings of the self-tangency and of the triple point strata of the discriminant.
Let $I$ be the set of moments when we crossed the triple point stratum.  
Let $\{\sigma_i\}_{i\in I}$ be the signs of the corresponding 
new born vanishing triangles and $\{[s_i]\}_{i\in I}$ be the $T$-equivalence 
classes represented by the corresponding generic 
curves with a triple point.
Put $\Delta_{\overline{\St}}(p)=\sum_{i\in I}\sigma_i\psi([s_i])$ and  
$\overline{\St}(\xi')=\overline{\St}(\xi)+\Delta_{\overline{\St}}(p)$.
To prove the Theorem it is sufficient to show, that $\overline{\St}(\xi')$ 
does not depend on the generic path $p$, we used to define it. 
The last statement follows from Lemma~\ref{onlyhomotop} and
Lemma~\ref{zerojump}.
Thus, we have proved Theorem~\ref{pfbarST} modulo these two lemmas.
\qed

\begin{lem}[Cf. V.~Arnold~\cite{Arnold}]\label{onlyhomotop}
Let $p$ be a generic path in $\mathcal F$, which connects $\xi$ to itself.
Then $\Delta_{\overline{\St}}(p)$ depends only on the class in $\pi_1(\mathcal
F, \xi)$ represented by $p$.
\end{lem}

\begin{lem}\label{zerojump}
If statement \textrm{II} of Theorem~\ref{barST} is true, then 
for every element of $\pi_1(\mathcal F, \xi)$ there exists a generic loop $q$
in 
$\mathcal F$, representing this element, such that $\Delta_{\overline{\St}}(q)=0.$
\end{lem}

\subsection{Proof of Lemma~\ref{onlyhomotop}}
It is sufficient to show that, if we go around any stratum of codimension two 
along a small generic loop $r$ (not necessarily starting at $\xi$), then 
$\Delta_{\overline{\St}}(r)=0$. 
The only strata of codimension two in the bifurcation
diagram of which triple points are present are: 
a) two distinct triple points, b) triple point and distinct self-tangency
point, c) triple point at which two
branches are tangent (of order one) and d) quadruple point (at which every
two branches are transverse). (All the codimension two singularities and
bifurcation diagrams for them were described by V.~Arnold~\cite{Arnold}.)

If $r$ is a small loop which goes around the stratum of two
distinct triple points, then in $\Delta_{\overline{\St}}(r)$ 
we have each of
the two $T$-equivalence classes twice, once with the plus sign of the
newborn vanishing triangle, once with the minus. 
Hence $\Delta_{\overline{\St}}(r)=0$.

If $r$ is a small loop which goes around the stratum of one triple and one
self-tangency point, then the two $T$-equivalence classes participating 
in $\Delta_{\overline{\St}}(r)$ are equal and the signs with which they
participate are opposite.
Hence $\Delta_{\overline{\St}}(r)=0$.

\begin{figure}[htbp]
 \begin{center}
  \epsfxsize 8.5cm
  \hepsffile{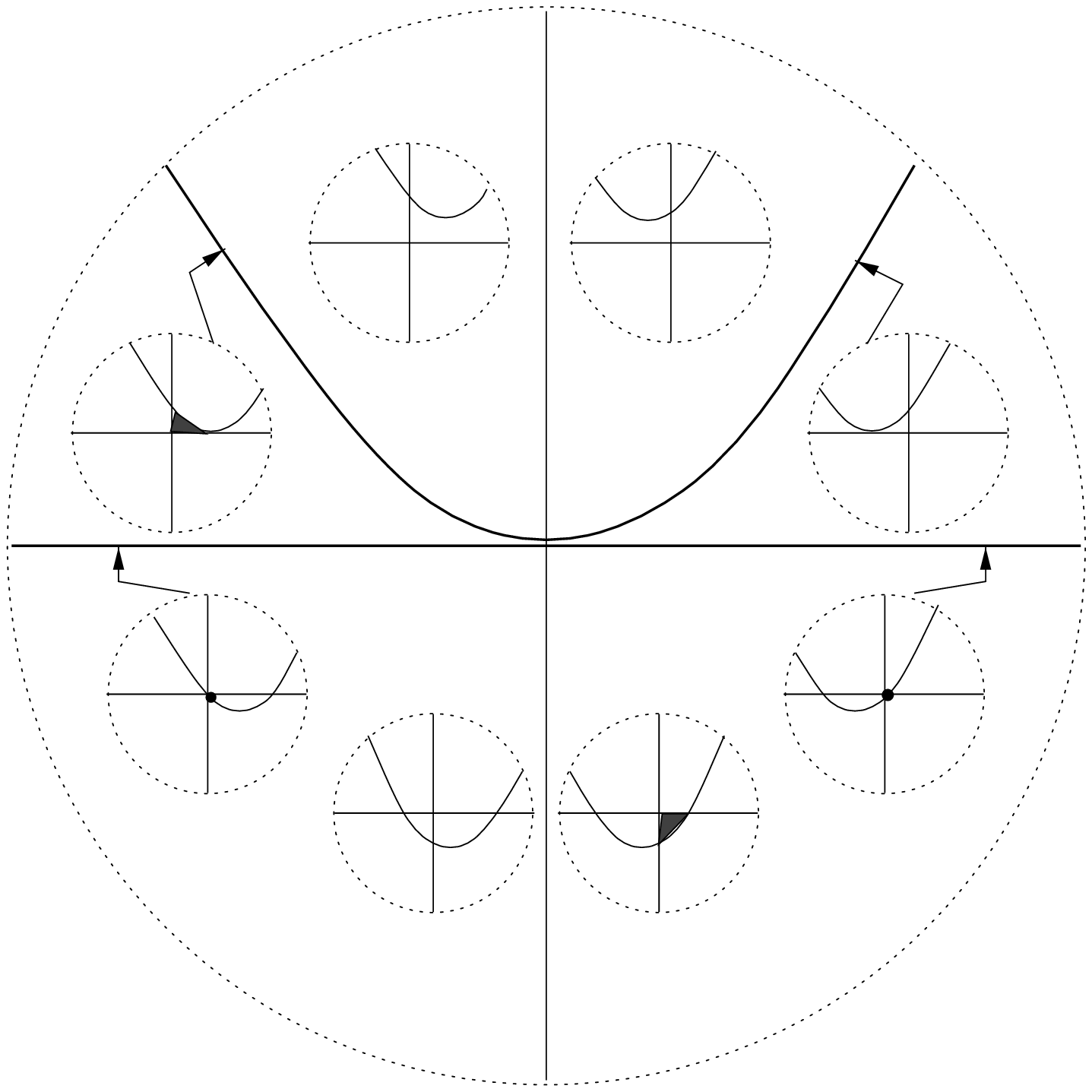}
 \end{center}
\caption{}\label{vers2.fig}
\end{figure}

Let $r$ be a small loop which goes around the stratum of a 
triple point with two
tangent branches. We can assume, that it corresponds to a loop on
Figure~\ref{vers2.fig} directed clockwise. (The colored triangles are the
newborn vanishing triangles.)
As we can see from Figure~\ref{vers2.fig} there are just two terms 
in $\Delta_{\overline{\St}}(r)$. 
It is clear, that the $T$-equivalence classes in them
coincide. A direct check shows that the signs of the two terms are
opposite. (Note, that if they are not always opposite, then
Arnold's $\St$ invariant is not well defined.)

\begin{figure}[htbp]
 \begin{center}
  \epsfxsize 8.5cm
  \hepsffile{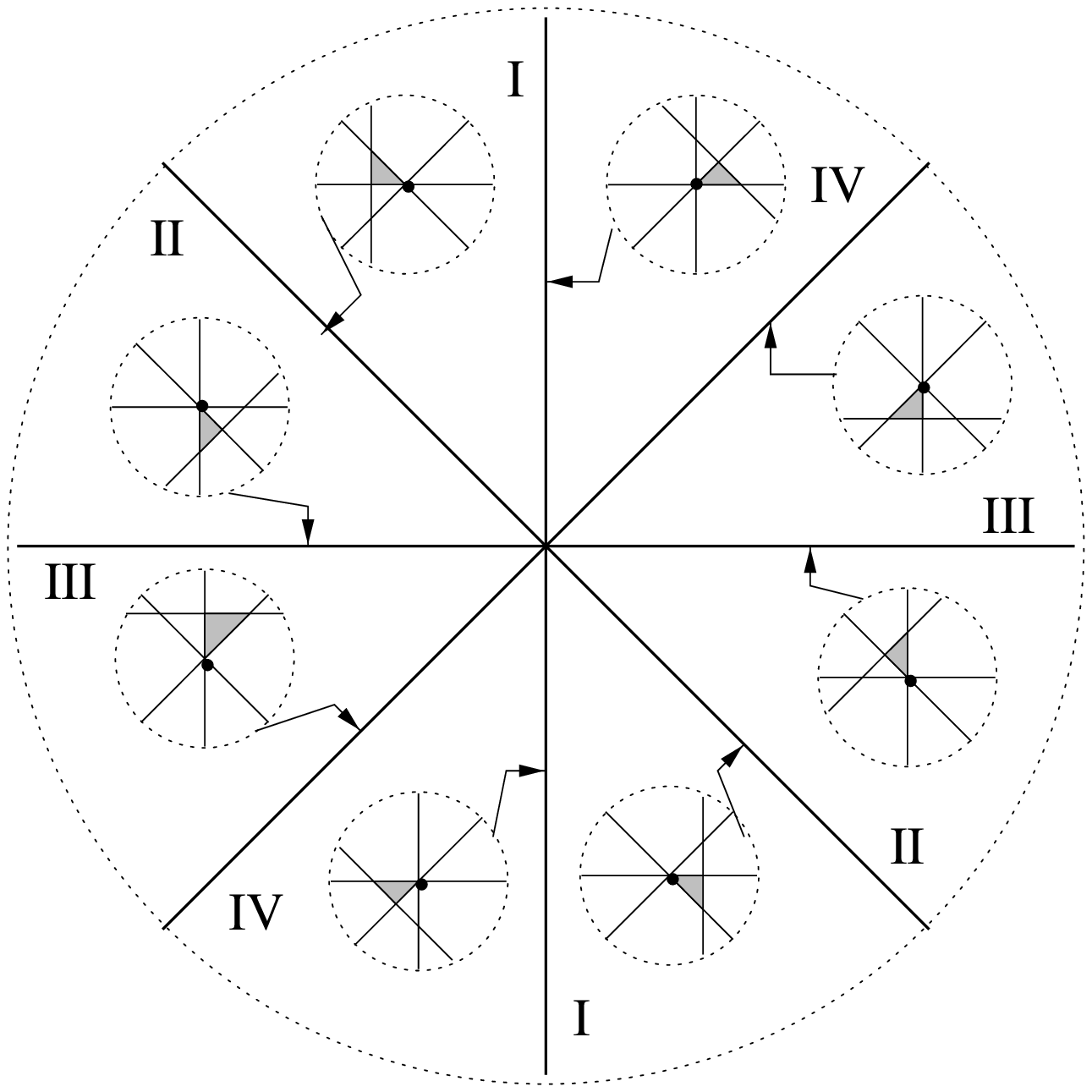}
 \end{center}
\caption{}\label{vers1.fig}
\end{figure}

Finally, let $r$ be a small loop, which goes around the stratum of a quadruple
point (at which every two branches are transverse). 
We can assume, that it corresponds to a loop in 
Figure~\ref{vers1.fig} directed counter clockwise. 
There are eight terms in $\Delta_{\overline{\St}}(r)$. We split them
into pairs I, II, III, IV, as it is shown in Figure~\ref{vers1.fig}. One can 
see, that 
the $T$-equivalence classes of the two curves in each pair are the same. 
For each branch the sign of the colored triangle 
is equal to the sign of the triangle, which died under the
triple point stratum crossing shown on the next (in the counterclockwise direction)
branch. 
The sign of the dying vanishing triangle 
is minus the sign of the newborn vanishing triangle.
Finally, one can see that the signs
of the colored triangles inside each pair are opposite. Thus, all these eight
terms cancel out.

This finishes the Proof of Lemma~\ref{onlyhomotop}.
\qed

\subsection{Proof of Lemma~\ref{zerojump}}
\begin{emf}\label{thoughtful}
In $\Z$ there are no elements of finite order. Thus,  
if $m\neq 0$, then $\Delta_{\overline{\St}}(q)\neq 0
\Leftrightarrow m\Delta_{\overline{\St}}(q)=
\Delta_{\overline{\St}}(q^m)\neq 0$. 
Hence, to prove Lemma~\ref{zerojump} 
it is sufficient to show that $\Delta_{\overline{\St}}(q^m)=0$ for a
certain power $m\neq 0$ of $q\in \pi_1(\mathcal F,\xi)$.
\end{emf}

\begin{prop}\label{commute} 
Let $F$ be a surface, $STF$ be its spherical tangent bundle 
and $p\in STF$ be a point.
Let $f\in \pi_1(STF,p)$ be
the class of an oriented (in some way) 
fiber of the $S^1$-fibration $\pr: STF\rightarrow F$. 

If $\alpha\in \pi_1(STF,p)$ is a loop projecting to an orientation
preserving loop on $F$, then 
\begin{equation}\label{commute1} 
\alpha f=f\alpha.
\end{equation}

If $\alpha\in \pi_1(STF,p)$ is a loop projecting to an orientation
reversing loop on $F$, then 
\begin{equation}\label{commute2}
\alpha f=f^{-1}\alpha.
\end{equation}
\end{prop}

The proof of this Proposition is straightforward.

\begin{emf}\label{H-principle}{\em Parametric $h$-principle.\/}
The parametric
$h$-principle, see~\cite{Gromov} page 16, implies
that $\mathcal F$ is weak homotopy equivalent to the space $\Omega STF$ 
of free loops in $STF$. 
The corresponding mapping $h:\mathcal F\rightarrow \Omega STF$ 
sends a curve $\xi\in \mathcal F$ 
to a loop $\vec \xi \in \Omega STF$ by  mapping a point $y\in S^1$ 
to a point in $STF$, to which points the 
velocity vector of $\xi$ at $\xi(y)$. 

Fix a point $a$ on $S^1$ (which parameterizes the curves).
Let $q$ be a loop in $\mathcal F$ starting at $\xi$. At any moment of time $q(t)$
is a curve, which can be lifted to a loop in $STF$. Thus, $q$ gives
rise to the mapping $q_h:S^1\times S^1\rightarrow STF$ (the lifting of $q$
by $h$). (In the product $S^1\times S^1$ the first copy of $S^1$ corresponds
to the parameterization of a curve and the second to the parameterization of
the loop $q$.) The mapping $q_h$ restricted to $a\times S^1$ gives rise to
the loop $t_a(q)$ in $STF$. (It is a trajectory of the lifting of $a$.)
One can check, that the mapping $t_a:\pi_1(\mathcal F, \xi)\rightarrow \pi_1(STF,
\vec \xi(a))$ is a homomorphism.

Note, that if $q\in \pi_1(\mathcal F, \xi)$ 
is the sliding of a kink (see~\ref{gamma1}) along a curve $\xi$ representing
an orientation preserving loop on $F$, 
then the velocity vector of $\xi$ at
$\xi(a)$ is rotated by $2\pi$ under this sliding. 
Thus, $t_a(q)\in\pi_1(STF, \vec\xi(a))$ is equal to 
$f$, the homotopy class of the fiber of the
$S^1$-fibration $\pr:STF\rightarrow F$.

One can check, that if $q\in \pi_1(\mathcal F, \xi)$ 
is the sliding of a kink along a curve $\xi$ representing
an orientation reversing loop on $F$, 
then $t_a(q)=1\in\pi_1(STF,\vec\xi(a))$.
(In this case the kink has to slide twice along $\xi$, before it
returns to its original position and 
the total angle of rotation of the velocity
vector of $\xi$ at $\xi(a)$ appears to be zero.) 
\end{emf}

\begin{prop}[Cf.~V.~Hansen~\cite{Hansen}]\label{G1}
The group $\pi_1(\Omega STF, \lambda)$ is isomorphic to $Z(\lambda),$ 
the centralizer of the element $\lambda\in \pi_1(STF, \lambda(a))$.
\end{prop}

\begin{emf}{\em Proof of Proposition~\ref{G1}.\/}
Let $p:\Omega STF\rightarrow
STF$ be the mapping, which sends $\omega\in\Omega STF$ to $\omega (a)\in
STF$. (One can check, that this $p$ is a Serre fibration, with
the fiber of it isomorphic to the space of loops
based at the corresponding point.)

A Proposition proved by V.L.~Hansen~\cite{Hansen} says that: if $X$ is a
topological space with $\pi_2(X)=0$, then
$\pi_1(\Omega X, \lambda)=Z(\lambda)<\pi_1(X, \lambda (a))$. (Here $\Omega X$
is the space of free loops in $X$ and $\lambda$ is an element of $\Omega X$.)
One can check that $\pi_2(STF)=0$ for any surface $F$. Thus, we get
that $\pi_1(\Omega STF, \lambda)$ is isomorphic to $Z(\lambda)<\pi_1(STF,
\lambda(a))$. From the proof of the Hansen Proposition it follows that the
isomorphism is induced by $p_*$. \qed
\end{emf}
 
The following statement is an immediate consequence of Proposition~\ref{G1} and
the $h$-principle (see Section~\ref{H-principle}).

\begin{cor}\label{center}
Let $F$ be a surface and $\xi$ be a curve on $F$, then $\pi_1(\mathcal F, \xi)$
is isomorphic to $Z(\vec\xi)$, the centralizer of $\vec\xi\in\pi_1(STF,
\vec\xi(a))$. The isomorphism is given by $t_a:\pi_1(\mathcal F,
\xi)\rightarrow Z(\vec\xi)$, which sends $q\in \pi_1(\mathcal F,
\xi)$ to $t_a(q)$. (See Section~\ref{H-principle}.)
\end{cor}

\begin{prop}\label{Preissman}
Let $F\neq S^2, T^2\text{ (torus), } \R P^2, K\text{ (Klein bottle)}$ 
be a surface (not necessarily compact or orientable)  
and $G'$ be a nontrivial commutative subgroup of $\pi_1(F)$. 
Then $G'$ is infinite cyclic 
and there exists a
unique maximal infinite cyclic $G<\pi_1(F)$, such that $G'<G$.
\end{prop}

\begin{emf}{\em Proof of Proposition~\ref{Preissman}.\/}
It is well known, that any closed $F$, other than $S^2, T^2, \R P^2, K$,
admits a hyperbolic metric of a constant negative curvature. 
(It is
induced from the
universal covering of $F$ by the hyperbolic plane $H$.) 
The Theorem by A.~Preissman (see~\cite{Docarmo}) 
says, that if $M$ is a compact Riemannian manifold with a negative curvature,
then any nontrivial Abelian subgroup $G'<\pi_1(M)$ is isomorphic to $\Z$.
Thus, if $F\neq S^2, T^2, \R P^2, K$ is closed, then any nontrivial 
commutative $G'<\pi_1(F)$ is infinite cyclic.

The proof of the Preissman Theorem given in~\cite{Docarmo} is based on
the fact, that if $\alpha,\beta\in\pi_1(M)$ are nontrivial commuting
elements, then there exists a
geodesic in $\bar M$ (the universal covering of $M$) which is mapped to
itself
under the action of these elements considered as deck transformations on
$\bar M$. Moreover, these transformations restricted to the geodesic act as
translations. This implies, that if $F\neq S^2, T^2, \R P^2, K$ 
is a closed surface, then there exists a unique maximal infinite cyclic
$G<\pi_1(F)$, such that $G'<G$. This gives the proof of
Proposition~\ref{Preissman} for
closed $F$. 

If $F$ is not closed then the statement of the Proposition is also
true, because in this case $F$ is homotopy equivalent to a bouquet of
circles.
\qed
\end{emf}

\centerline{}

We first prove Lemma~\ref{zerojump} 
for $F\neq S^2, \R P^2, T^2, K$ and then separately 
for the cases
$F=S^2, \R P^2, T^2, K$.

\begin{emf}\label{other}{\em Case $F\neq S^2, T^2, \R P^2, K$\/}
Corollary~\ref{center} says, that 
$\pi_1(\mathcal F, \xi)=Z(\vec\xi)<\pi_1(STF, \vec\xi(a))$. The corresponding
isomorphism (see Section~\ref{H-principle}) 
maps $q\in\pi_1(\mathcal F, \xi)$ to $t_a(q)\in\pi_1(STF, \xi(a))$.

Thus, for any  $q\in\pi_1(\mathcal F, \xi)$ 
the elements $t_a(q)$ and $\vec \xi$ commute in $\pi_1(STF,
\vec \xi(a))$. Hence, $\xi=\pr_*(\vec \xi)$ 
commutes with $\pr_*(t_a(q))$ in $\pi_1(F,
\xi(a))$. Proposition~\ref{Preissman} implies, that 
there is exists an infinite cyclic subgroup of 
$\pi_1(F,\xi(a))$ generated by some $g\in
\pi_1(F,\xi(a))$, which contains both of these loops.  
Then there exist $m,n \in \Z$ such that $\xi=g^m$ and
$\pr_*(t_a(q))=g^n$.

Consider a curve $l$ (direct tangent to $\xi$ at $\xi(a)$) which represents $g$. 
We can lift it to an element $\vec g\in\pi_1(STF,
\vec \xi(a))$. 

The kernel of $\pr_*$ is generated by $f$, the class of an oriented fiber.
Using~\eqref{commute1} and~\eqref{commute2} one can interchange $f$ 
with the other elements of $\pi_1(STF,\vec \xi(a))$. We get that 
$\vec \xi=\vec g^mf^k$ and 
$t_a(q)=\vec g^nf^l$, for some $k,l\in \Z$. We prove Lemma~\ref{zerojump} separately for
cases $m\neq 0$
and $m=0$ in Section~\ref{nontrivial} and Section~\ref{trivial}, respectively.
(These two cases correspond to $\xi \neq 1\in\pi_1(F, \xi(a))$ and 
to $\xi=1\in\pi_1(F, \xi(a))$, respectively.) 
\end{emf} 

\begin{emf}\label{nontrivial}{\em Case $m\neq 0$.\/}
To prove Lemma~\ref{zerojump} it is sufficient to show, that 
$\Delta_{\overline{\St}}(q^m)\neq 0$
(see~\ref{thoughtful}). 

One can show that $t_a(q^m)=\vec\xi^nf^j$ for some $j\in \Z$. For $g$
which is an orientation preserving loop this follows from the
following calculation (which uses~\eqref{commute1}):  
\begin{equation}
t_a(q^m)=(t_a(q))^m=\bigl ( \vec g^n f^l\bigr )^m=
\bigl( \vec g^mf^k \bigr)^nf^{lm-nk}=\vec \xi^nf^{lm-nk}.
\end{equation}
For $g$, which is an orientation reversing loop
on $F$, this follows from the similar calculation (which
uses~\eqref{commute2}).
The fact that $t_a(q^m)$ should commute with $\vec \xi$ (since it is the 
$m$-th power of $t_a(q)\in Z(\vec\xi)$) 
and the identity~\eqref{commute2} imply,
that $t_a(q^m)=\vec \xi^n$, 
provided that $\xi$ represents an orientation reversing loop on $F$. 

Let $\gamma_1$ be the sliding of a kink along $\xi$ (see~\ref{gamma1}). 
If $\xi$ represents an orientation preserving loop on $F$, 
then the velocity vector of $\xi$ at $\xi(a)$ 
is rotated by $2\pi$ under $\gamma_1$. 
Thus, 
$t_a(\gamma_1)=f$.
Hence, the loop $\alpha\in\pi_1(\mathcal F, \xi)$ for which 
$t_a(\alpha)=t_a(q^m)$ is: $n$ times sliding of $\xi$
along itself according to the orientation, composed with
$\gamma_1^{j}$.

As it was said above, if $\xi$ represents an orientation
reversing loop, then $t_a(q^m)=\vec\xi^n$. 
Hence, the loop $\alpha\in\pi_1(\mathcal F, \xi)$, for which 
$t_a(\alpha)=t_a(q^m)$ is: $n$ times sliding of $\xi$
along itself. 
(In~\ref{H-principle} it was shown that $\gamma_1=1\in \pi_1(\mathcal F, \xi)$ 
for $\xi$ representing an orientation reversing loop on $F$.)

No triple points appear
during the sliding of $\xi$ along itself. 
The inputs of the triple point stratum crossings which
occur under $\gamma_1$ cancel out, by the assumption of the Lemma. 
Hence, $\Delta_{\overline{\St}}(q^m)=0$.

Thus, we have proved (see~\ref{thoughtful}) Lemma~\ref{zerojump} 
for $F\neq S^2, \R P^2, T^2, K$ and $m\neq 0$. 
\end{emf}

\begin{emf}\label{trivial}{\em Case $m=0$.\/}
If $m=0$, then $\xi$ represents $1\in\pi_1(F,\xi(a))$. 
For any $q\in\pi_1(\mathcal F, \xi)$ the projection of
$t_a(q^2)\subset STF$ to $F$ is
an orientation preserving loop on $F$.
A straightforward check shows that for any $q\in\pi_1(\mathcal F, \xi)$ the
element $q^2$ can be obtained by a composition of 
$\gamma_1^{\pm 1}$
(see~\ref{gamma1}) and loops obtained by the following construction.
  
Push $\xi$ into a small disc by a generic regular homotopy $r$. 
Slide
this small disc along some orientation preserving curve in $F$ 
and return $\xi$ to its original
shape along $r^{-1}$. 

Clearly, the inputs of $r$ and $r^{-1}$ into 
$\Delta_{\overline{\St}}$ cancel out and no triple point stratum
crossings happen, when
we slide a small disc along a path in $F$. Thus, loops obtained by this
construction do not give any input to $\Delta_{\overline{\St}}$.
By the assumption of the Lemma $\Delta_{\overline{\St}}(\gamma_1)=0.$ 

This implies that $\Delta_{\overline{\St}}(q^2)=0$
for any $q\in \pi_1(\mathcal F, \xi)$, and
we have proved (see~\ref{thoughtful}) Lemma~\ref{zerojump}  
for $F\neq S^2, \R P^2, T^2, K$. 
\end{emf}

\begin{emf}\label{S^2}{\em Case $F=S^2$.\/}
One checks that $\pi_1(STS^2)=\Z_2$. (Note that $STS^2=\R P^3$.) 
Corollary~\ref{center} implies that 
$\pi_1(\mathcal F, \xi)=\Z_2$ for $F=S^2$.
Thus, $\Delta{\overline{\St}}(q^2)=\Delta{\overline{\St}}(1)=0$. (Here $1$ is
a trivial loop in $\mathcal F$.) This finishes (see~\ref{thoughtful}) the proof 
of Lemma~\ref{zerojump} for $F=S^2$. 
\end{emf}
 
\begin{emf}\label{T^2}{\em Case $F=T^2$.\/}
Using identity~\eqref{commute1} we get, that
$\pi_1(STT^2)=\Z\oplus\Z\oplus\Z$. 
Corollary~\ref{center} implies that $\pi_1(\mathcal F, \xi)
=\pi_1(STT^2)=\Z\oplus\Z\oplus\Z$.
The generators of this group are:

1) The loop $\gamma_1$, which is the sliding of a kink along 
$\xi$ (see~\ref{gamma1}). 

2) The loops $\gamma_3$ and $\gamma_4$, which are  
slidings of $\xi$ along the 
unit vector fields parallel to the meridian and longitude
of $T^2$, respectively.

By the assumption of the Lemma 
$\Delta_{\overline{\St}}(\gamma_1)=0$. Since no discriminant 
crossings occur during $\gamma_3$ and $\gamma_4$ we get, that 
$\Delta_{\overline{\St}}(\gamma_3)=\Delta_{\overline{\St}}(\gamma_4)=0$.
This finishes the proof of Lemma~\ref{zerojump} for $F=T^2$.
\end{emf}

\begin{emf}\label{RP^2}{\em Case $F=\R P^2$.\/}
One checks that $\pi_1(ST \R P^2)=\Z_4$. 
Corollary~\ref{center} implies that 
$\pi_1(\mathcal F, \xi)=\Z_4$ for $F=\R P^2$.
Thus, $\Delta{\overline{\St}}(q^4)=\Delta{\overline{\St}}(1)=0$. (Here $1$ is
a trivial loop in $\mathcal F$.) This finishes (see~\ref{thoughtful}) the proof 
Lemma~\ref{zerojump} for $F=\R P^2$. 

\end{emf}

\begin{rem}
As it follows from~\ref{RP^2} the condition
$\Delta{\overline{\St}}(\gamma_1)=0$ is automatically satisfied in the case
of curves on $\R P^2$.
\end{rem}

\begin{emf}\label{K}{\em Case $F=K$.\/}
Corollary~\ref{center} says, that 
$\pi_1(\mathcal F, \xi)$ is isomorphic to $Z(\vec \xi)<\pi_1(STK, \xi(a))$.

Consider $K$ as a quotient of a rectangle modulo the identification on its
sides, which is shown in Figure~\ref{klein.fig}.
We can assume that $\xi(a)$ coincides
with the image of a corner of the rectangle and that $\xi$ is direct
tangent to the
curve $c$ at $\xi(a)$. Let $g$ and $h$ be the
curves such that: $\vec \xi(a)=\vec g(a)=\vec h(a)$, $g=c\in \pi_1(K,
\xi(a))$ and
$h=d\in \pi_1(K, \xi(a))$. (Here $c$ and $d$ are the elements of $\pi_1(K)$
realized by the sides of the rectangle used to construct $K$,
see Figure~\ref{klein.fig}.)
Let $f$ be the
class of an oriented fiber of the fibration $\pr:STK\rightarrow K$.
One can show that:
\begin{equation}\label{present}
\pi_1(STK, \vec \xi(a))=\bigl \{\vec g, \vec h, f\big | \vec h\vec g^{\pm 1}
=\vec g^{\mp 1}\vec h, \text{ }\vec hf^{\pm
1}=f^{\mp 1}\vec h, \text{ }\vec g f=f\vec g\bigr \}.
\end{equation}

The second and the third relations in this presentation follow
from~\eqref{commute1} and~\eqref{commute2}.
To get the first relation one notes that                   
the identity $d c^{\pm 1}=c^{\mp 1} d\in\pi_1(K, \xi(a))$ implies
$\vec h \vec g^{\pm 1}=\vec g^{\mp 1} \vec h f^k$ for some $k\in \Z$.
But $\vec h^2$ commutes with $\vec g$, since they can be lifted to $STT^2$,
the fundamental group of which is Abelian. Hence, $k=0$.

Using relations~\eqref{present}                           
one can calculate $Z(\vec \xi)=\pi_1(\mathcal F, \xi)$.        
(Note that these relations allow one to present any        
element of $\pi_1(STK, \vec \xi(a))$ as                    
$\vec g^k\vec h^l f ^m$, for some $k,l,m\in \Z$.)

This group appears to be:

{\bf a)} The whole group $\pi_1(STK,\vec\xi(a))$, provided that $\vec
\xi=\vec h^{2l}$ for
some $l\in \Z$.

{\bf b)} An isomorphic to $\Z\oplus\Z\oplus\Z$
subgroup of $\pi_1(STK, \vec\xi(a))$, provided that $\vec\xi=\vec g^k\vec
h^{2l}f^m$ for some $k,l,m\in\Z$, such that $k\neq 0$ or $m\neq 0$.
This subgroup is generated by $\{\vec g, \vec h^2,f\}$.

{\bf c)} An isomorphic to $\Z\oplus\Z$ subgroup of $\pi_1(STK, \vec\xi(a))$,
provided that $\vec \xi=\vec g^k\vec h^{2l+1} f^m$ for some $k,l,m\in\Z$.
This subgroup is generated by $\{\vec \xi, \vec h^2\}$.

A straightforward check (which uses~\eqref{present}) 
shows that:

{\bf a)} If $\xi$ represents an orientation preserving loop on $K$, then a certain
degree of any loop $\gamma\in \pi_1(\mathcal F, \xi)$ can be expressed
as a product of $\gamma_1$ (see~\ref{gamma1}), 
$\gamma_2$ (see~\ref{gamma2}), $\gamma_3$, described below,  and their
inverses.

{\bf b)} If $\xi$ represents an orientation reversing loop on $K$, then a certain
degree of any loop $\gamma\in \pi_1(\mathcal F, \xi)$ can be expressed
as a product of $\gamma_3, \gamma_4$, described below, and their inverses.

Consider a loop $\beta$ 
in the space of all the autodiffeomorphisms of
$K$, which is the sliding of $K$ along the unit vector field parallel 
to the curve $d$ on $K$. 
(Note that $K$ has to slide twice along itself under this loop before all
points of $K$ come  
to the original position.) The loop $\gamma_3$ is the 
sliding of $\xi$ induced by $\beta$.

The loop $\gamma_4$ is the sliding of $\xi$ along itself.

No triple point stratum crossings occur under $\gamma_3$ and
$\gamma_4$. By the assumption of the Lemma
$\Delta_{\overline{\St}}(\gamma_1)=0$ and 
$\Delta_{\overline{\St}}(\gamma_2)=0$ (when $\gamma_2$ is well defined).

Thus, $\Delta_{\overline{\St}}(\gamma_i)=0$ ($i\in\{1,2,3,4\})$ 
and we have proved (see~\ref{thoughtful}) Lemma~\ref{zerojump} in the case
of $F=K$.

\begin{rem}
One can check that for the curve on the Klein bottle shown in
Figure~\ref{klein3.fig} the equations $\Delta_{\overline{\St}}(\gamma_1)=0$
and $\Delta_{\overline{\St}}(\gamma_2)=0$ are independent. This means that
both of the corresponding conditions are needed for the integrability of
$\psi$.

\begin{figure}[htbp]
 \begin{center}
  \epsfxsize 4cm
  \hepsffile{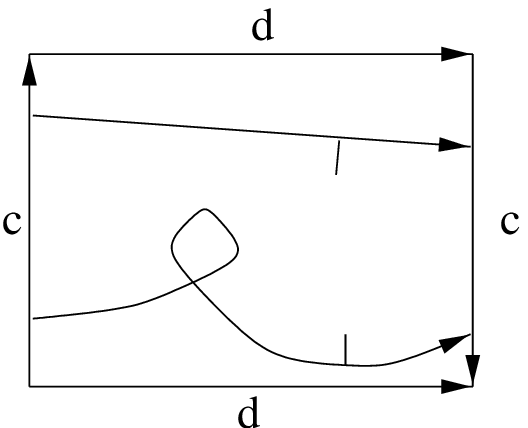}
 \end{center}
\caption{}\label{klein3.fig}
\end{figure}

\end{rem}

\end{emf}

This finishes the proof of Lemma~\ref{zerojump} for all the cases.
\qed

\subsection{Proof of Theorem~\ref{discrST}}\label{pfdiscrST}
Let $\mathcal S$ be the space of all the smooth mappings (not necessarily
immersions) from $S^1$ to $F$. 
Consider the subspace $\mathcal N$ of $S^1(3)\times \mathcal S$ consisting of 
$t\times f$, such that $f$ maps the three points from $t$ to one point on 
$F$. Clearly $\mathcal M$ is a subspace of $\mathcal N$.

Let $s_1, s_2\in \mathcal C$ be two generic curves with a triple point, such
that $[s_1]=[s_2]$. 
Let $m_1, m_2\in \mathcal M$ be the elements 
corresponding to $s_1$ and $s_2$.
To prove the Theorem we need to show that $m_1$ and $m_2$ belong
to the same path connected component of $\mathcal M$. 

\begin{figure}[htbp]
 \begin{center}
  \epsfxsize 9cm
  \hepsffile{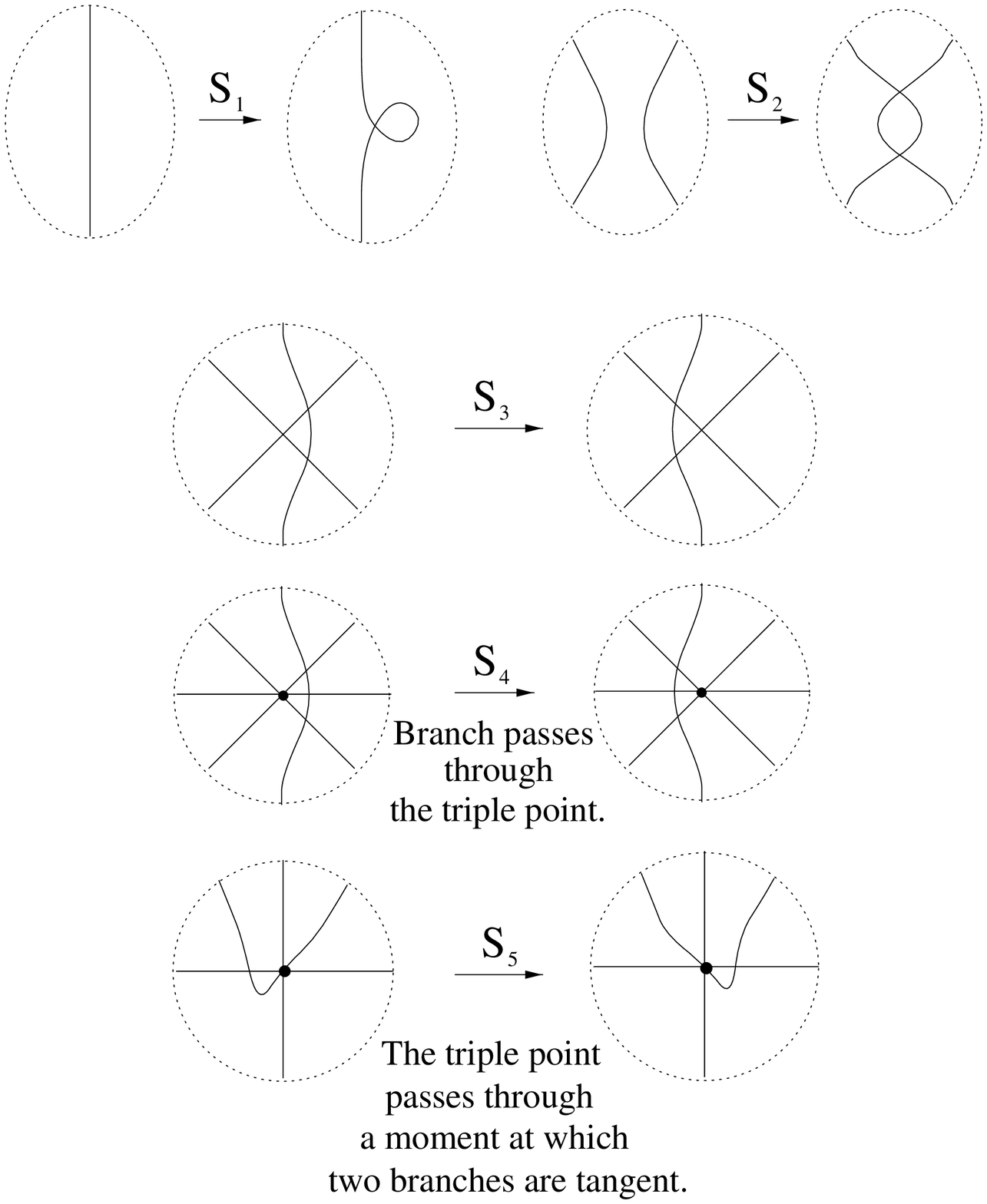}
 \end{center}
\caption{}\label{moves1.fig}
\end{figure}

Since $[s_1]=[s_2]$, we see that $m_1$  can be transformed to $m_2$ (in
$\mathcal N$) by the sequence of moves $S_1, S_2, S_3, S_4, S_5$ (shown in
Figure~\ref{moves1.fig}) and their inverses (and a continuous 
change of the parameterization). 
Note, that $S_1^{\pm}$ is 
the only move in
this list which happens not in $\mathcal M$. 

We can imitate the $S_1$-move
staying in $\mathcal M$ by creating two opposite kinks (see 
Figure~\ref{twokink.fig}) and then making one of these kinks very small.
(Similarly, we can imitate the $S_1^{-1}$-move staying in $\mathcal M$.)

We use $S_2, S_3, S_4, S_5$ and their inverses and the imitations of
$S_1^{\pm 1}$ to deform $m_1$ inside $\mathcal M$ so that
it looks exactly as $m_2$, except some number of small extra kinks
located on the three loops of $m_1$ adjacent to the triple point.

As it is shown below, we can create two opposite extra kinks, 
the first on one of the three loops of $m_1$, the second on another.

\begin{figure}[htbp]
 \begin{center}
  \epsfxsize 12cm
  \hepsffile{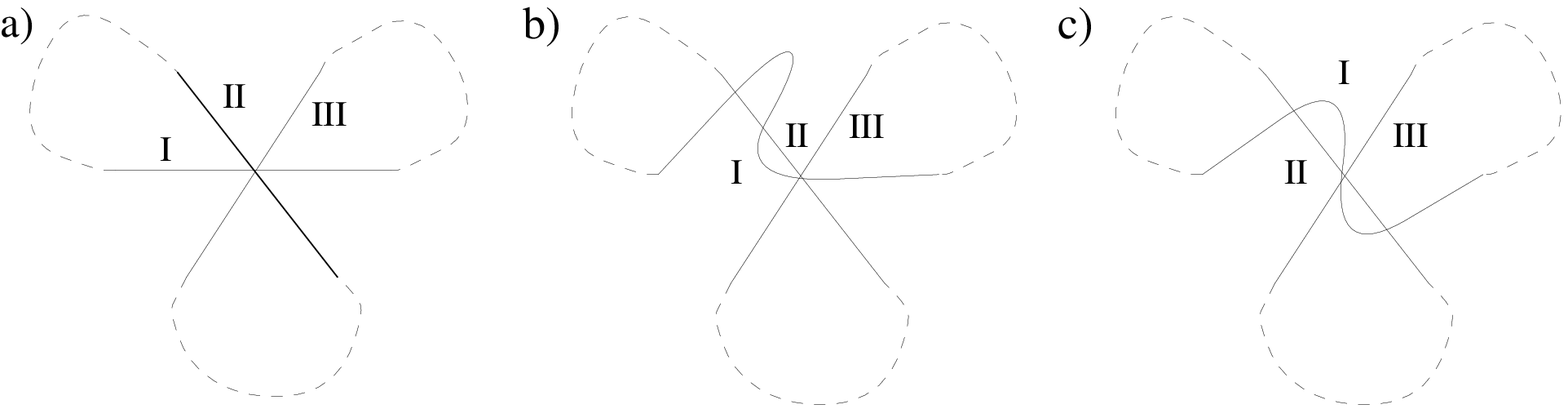}
 \end{center}
\caption{}\label{exchange.fig}
\end{figure}

The order, in which a small loop going around the triple point 
crosses the three branches of $m_1$ passing through the triple 
point, induces a
cyclic order on the branches. 
We use $S_2$ to deform Figure~\ref{exchange.fig}a to
Figure~\ref{exchange.fig}b, then we use $S_5$ to deform it to
Figure~\ref{exchange.fig}c. Note, that under this procedure the branches
$\textrm{I}$ and $\textrm{II}$ get interchanged in the cyclic order. 
Then in a
similar way we interchange the branches in the pairs $\{\textrm{I}, \textrm{III}\}$; 
$\{\textrm{I}, \textrm{II}\}$ and
$\{\textrm{I}, \textrm{III}\}$. 
After this the local picture around the triple point is the same as before.
One can check, 
that what happened with $m_1$ globally is equivalent to the addition of two
opposite kinks, the first to one branch of $m_1$, the second to another.  

It is clear, that 
using this procedure and the cancelation of two opposite kinks (see
Figure~\ref{twokink.fig}) we
can concentrate all the extra
kinks on one of the three loops of $m_1$. 
Slide these extra kinks along the loop, so that they are
all concentrated on a small arc.
Cancel out all the pairs of
opposite extra kinks. Now all the small extra kinks are
pointing to one side of the loop.

The $h$-principle (see~\cite{Gromov} page 16) 
implies that the space of all the curves on $F$ is
weak homotopy equivalent to the space of all the free loops in $STF$ (the
spherical tangent bundle of $F$). 
The corresponding mapping $h$ sends a curve $\xi\in \mathcal F$ 
to a loop $\vec\xi\subset STF$ by  mapping a point $y\in S^1$ 
to a point in $STF$, to which points the 
velocity vector of $\xi$ at $\xi(y)$. Since by our assumption 
both $s_1$ and $s_2$ belong to
the same connected component of $\mathcal F$, we get that 
their liftings to loops in $STF$ are free homotopic.

Let $f$ be the homotopy class of the fiber of the $S^1$-fibration
$\pr:STF\rightarrow F$. 
An extra small kink corresponds under $h$ to a multiplication by 
$f^{\pm 1}$, depending on
the side of the loop the kink points to. Let $n$ be the number of small extra
kinks which are present on $m_1$.

Fix a point $a$ on $S^1$. We can assume that after the process described
above the curves $s_1$ and $s_2$ (corresponding to $m_1$ and $m_2$,
respectively) are direct tangent at the image of $a$.  
Now we can consider $s_1$ and $s_2$ as elements of $\pi_1(F, s_1(a))$ and 
the liftings $\vec s_1$ and $\vec s_2$ 
(see~\ref{H-principle}) as elements of $\pi_1(STF, \vec s_1(a))$. By the
initial assumption $s_1$ and $s_2$ belong to the same connected component 
of $\mathcal F$. 
The $h$-principle implies, that $\vec s_1$ is free homotopic to 
$\vec s_2$.  Hence, we get that for some element $\alpha\in\pi_1(STF, s_1(a))$ 
\begin{equation}\label{kinknumber1}
\vec s_1=\alpha \vec s_1 f^n \alpha^{-1}.
\end{equation} 

Consider the case of $F=S^2$. One checks that $\pi_1(STS^2)=\Z_2$ is
commutative and $f$ has order two in $\pi_1(STS^2)$. (Note that $STS^2=\R
P^3$.)
From~\eqref{kinknumber1} we get that $n$ (the number of extra kinks) is even.
We take one of the kinks and evert it, by expanding it
till it goes around $S^2$ and comes back as a kink pointing to the other
side of the loop. Then we cancel it out with one of the other extra kinks.
In order to deform $m_1$ to $m_2$ we perform this operation until
there are no extra kinks left.

In the case of $F=T^2$ the group $\pi_1(STT^2)=\Z\oplus\Z\oplus\Z$ is
commutative. From~\eqref{kinknumber1} we get that $f^n=1$. But 
$f\in\pi_1(STT^2)$ has infinite order, thus $n=0$ and there were no extra
kinks that survived the process. This means, that we have
constructed the desired path from $m_1$ to $m_2$.

For $F\neq S^2, T^2$ the element $f\in\pi_1(STF)$ has infinite order.
Combining  
identities~\eqref{kinknumber1} and~\eqref{commute1} 
(recall that $F$ was assumed to be orientable) we
get that $\vec s_1^{-1}\alpha^{-1}\vec s_1 \alpha=f^n$. 
Thus, the projections
of $\vec s_1$ and $\alpha$ commute in $\pi_1(F, s_1(a))$.
Proposition~\ref{Preissman} implies that these projections can be expressed
as powers of some $g\in\pi_1(F, s_1(a))$. 
Let $g_{s_1}$ be a curve representing this $g$, which is 
direct tangent to 
$s_1$ at $s_1(a)$.
The kernel of the homomorphism $\pr_*$ is generated by $f$. Using
identity~\eqref{commute1}  
we can present $\alpha$ as $\vec g_{s_1}^i f^j\in\pi_1(STF,
\vec s_1(a)$ and $s_1$ as $\vec g_{s_1}^k f^l\in \pi_1(STF, \vec s_1(a))$,
for some $i,j,k,l\in\Z$. 
This means (see~\eqref{commute1}) that 
$\alpha$ commutes with $s_1$ in $\pi_1(STF, \vec s_1(a))$. From the 
identity~\eqref{kinknumber1} we get that $f^n=1$. But $f$ has infinite order
in $\pi_1(STF)$. Hence, $n=0$ and
there were no extra kinks, that survived the process. 
This means that we have
constructed the desired path from $m_1$ to $m_2$.
\qed

\subsection{Proof of Theorem~\ref{barJ+}.}\label{pfbarJ+}
The proof of Theorem~\ref{barJ+} is analogous to the proof of
Theorem~\ref{barST}.

One can easily formulate and prove the corresponding versions of 
Lemma~\ref{onlyhomotop}
and Lemma~\ref{zerojump}.

\begin{figure}[ht]
 \begin{center}
  \epsfxsize 8.5cm
  \hepsffile{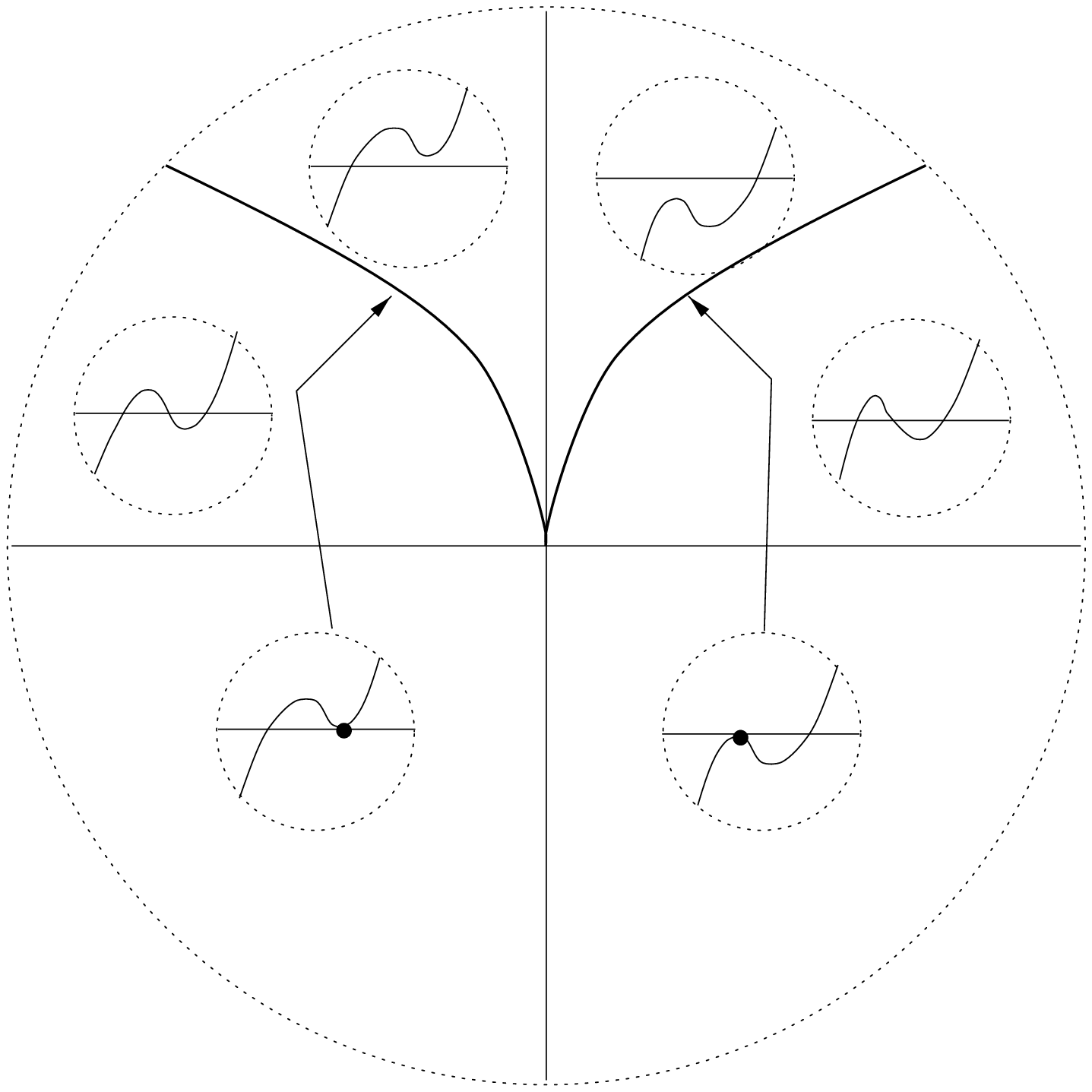}
 \end{center}
\caption{}\label{vers3.fig}
\end{figure}

The strata you have to go around, in order to prove the analogue of
Lemma~\ref{onlyhomotop}, are: a) two self-tangency points, b) self-tangency
point and distinct triple point,
c) triple point at
which exactly two branches are tangent (of order one), 
d) self-tangency point of order two. 
The bifurcation diagrams for the last two cases are shown in
Figure~\ref{vers2.fig} and Figure~\ref{vers3.fig}, respectively.
\qed

\subsection{Proof of Theorem~\ref{discrJ+}}\label{pfdiscrJ+}
The proof of this Theorem is analogous to the proof of
Theorem~\ref{discrST}. 

\begin{figure}[htbp]
 \begin{center}
  \epsfxsize 7cm
  \hepsffile{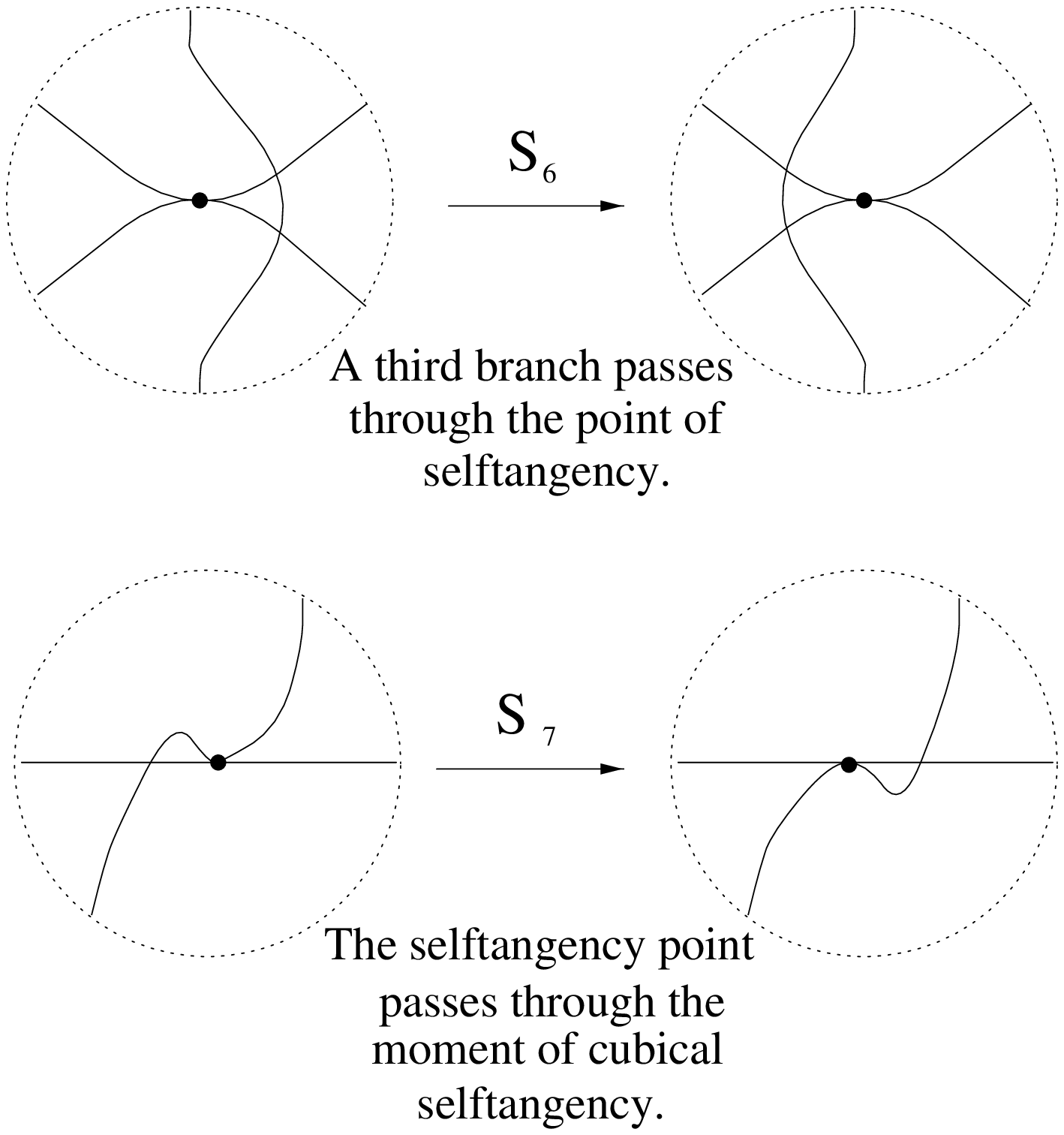}
 \end{center}
\caption{}\label{moves2.fig}
\end{figure}

Let $s_1,s_2\in \mathcal C$ be two generic curves with a point of direct
self-tangency, 
%(see~\ref{genericdirect}), 
such that $[s^+_1]=[s^+_2]$.
Let $m_1, m_2\in \mathcal M^+$ be the elements corresponding to $s_1$ and $s_2$,
respectively. Since $[s^+_1]=[s^+_2]$ we can choose 
the mappings of $B_2$ to $STF$ 
associated with $s_1$ and $s_2$ so that they are free homotopic. 
Hence, the
projections of them to $F$ are also free homotopic. One can show, that 
the projections of the two circles of $B_2$ can be assumed to be 
direct tangent (at the base point) under this homotopy. 
Clearly the only moves needed for this homotopy are $S_1, S_2,
S_3$ (see Figure~\ref{moves1.fig}), $S_6, S_7$ (see Figure~\ref{moves2.fig})
and their inverses. 

We use the imitation of $S_1$ (described in~\ref{pfdiscrST}) $S_3, S_4,
S_6, S_7$ and
their inverses (and a continuous change of the parameterization) 
to deform $m_1$ in the space $\mathcal M^+$ to an element, which
looks nearly as $m_2$, except a number of small extra kinks located on the two
loops of $m_1$ adjacent to the point of a direct self-tangency.

For each of the two loops we slide all the extra kinks, so that they are
located on a small arc of the loop. We cancel all the pairs of opposite
kinks by reversing the process shown in Figure~\ref{twokink.fig}. 
Now the kinks on each loop are pointing to the same side of it. 

We note,
that the $T^+$-equivalence class 
corresponding to $m_1$ did not change under all
these deformations. 
An extra kink located on a loop of $m_1$ corresponds (under the lifting 
of the loop to a loop in $STF$) to the multiplication
by the class of an oriented fiber of $\pr:STF\rightarrow F$. 
Similarly to~\ref{pfdiscrST} we get, that for $F\neq
S^2$, the number of extra kinks on each of the two loops of $m_1$ is 
zero. This means that we have constructed the desired 
path connecting $m_1$ to $m_2$. 
For $F=S^2$  we use the process 
described in Section~\ref{pfdiscrST} 
to cancel out all the extra kinks on each
of the two loops, and obtain a path connecting $m_1$ to $m_2$.

This finishes the Proof of Theorem~\ref{discrJ+}.
\qed

\centerline{}
\centerline{\bf Acknowledgments}
The main results contained in this paper formed part of my Ph.D.
thesis~\cite{Tchernov3}
(Uppsala University 1998). 
I am deeply grateful to Oleg Viro, Tobias Ekholm and Viktor Goryunov
for many enlightening discussions.

\end{document}